\documentclass{amsart} %{amsart}%{article}%{elsart}
\usepackage{graphicx}
\usepackage{amssymb, amsmath, sidecap}
\usepackage{float}
\usepackage{extarrows}
\usepackage{mathrsfs}
\usepackage{booktabs}
\usepackage{verbatim}
\usepackage{hyperref}
\usepackage[usenames,dvipsnames]{xcolor}

\renewcommand{\baselinestretch}{1}

\def\bt{\begin{thm}}
\def\et{\end{thm}}
\def\bl{\begin{lem}}
\def\el{\end{lem}}
\def\bd{\begin{defi}}
\def\ed{\end{defi}}
\def\bc{\begin{cor}}
\def\ec{\end{cor}}
\def\bp{\begin{proof}}
\def\ep{\end{proof}}
\def\br{\begin{rem}}
\def\er{\end{rem}}

\def\Forall{\text{ } \forall \:}
\def\d{\, \mathrm{d}}

\def\be{\begin{equation}}
\def\ee{\end{equation}}
\def\bes{\begin{equation*}}
\def\ees{\end{equation*}}
\def\bea{\begin{equation} \begin{aligned}}
\def\eea{\end{aligned} \end{equation}}
\def\beas{\begin{equation*} \begin{aligned}}
\def\eeas{\end{aligned} \end{equation*}}
\def\ba{\begin{align}}
\def\ea{\end{align}}
\def\bas{\begin{align*}}
\def\eas{\end{align*}}

\newtheorem{thm}{Theorem}[section]
\newtheorem{lem}{Lemma}[section]
\newtheorem{defi}{Definition}[section]

\newtheorem{prop}[thm]{Proposition}

\newtheorem{rem}{Remark}[section]
\newtheorem{cor}{Corollary}[section]

\numberwithin{equation}{section}
\numberwithin{figure}{section}

\begin{document}
\title{H\"{o}lder continuity for nonlinear elliptic problem in Musielak-Orlicz-Sobolev space \footnotemark[1]\footnotemark[2]}

\author[Wang]{Beibei Wang}
\address[Beibei Wang]{School of Mathematics and Statistics, Lanzhou
University, Lanzhou 730000, P. R. China} \email{wangbb15@lzu.edu.cn}

\author[Liu]{Duchao Liu}
\address[Duchao Liu]{School of Mathematics and Statistics, Lanzhou
University, Lanzhou 730000, P. R. China} \email{liuduchao@gmail.com,
Tel.: +8613893289235, fax: +8609318912481}

\author[Zhao]{Peihao Zhao}
\address[Peihao Zhao]{School of Mathematics and Statistics, Lanzhou
University, Lanzhou 730000, P. R. China} \email{zhaoph@lzu.edu.cn}

\footnotetext[1]{
In memory of Professor Wenyuan Chen.}

\footnotetext[2]{
Research supported by the National Natural Science Foundation of
China (NSFC 11501268 and NSFC 11471147).
%and the Fundamental Research Funds for the Central Universities (lzujbky-2016-101).
}

\begin{abstract}
Under appropriate assumptions on the $N(\Omega)$-fucntion, the De Giorgi process
is presented in the framework of Musielak-Orlicz-Sobolev space to prove the H\"{o}lder continuity of fully nonlinear elliptic problems. As the applications, the
H\"{o}lder continuity of the minimizers for a class of the energy
functionals in Musielak-Orlicz-Sobolev spaces is proved; and furthermore, the H\"{o}lder continuity of the weak solutions for a class of fully nonlinear elliptic equations is provided.
\end{abstract}
 \maketitle

%\tableofcontents

\section{Introduction}

Since Ladyzhenskaya and Ural'tseva developed the method pioneered by De Giorgi \cite{DeGiorgi57} and introduced the class $\mathscr{B}(\Omega,M,\gamma,\gamma_{1},\delta,1/q)$ (see \cite{Ladyzhenskaya}), through which the H\"{o}lder continuity of functions of this class can be proved, the class $\mathscr{B}(\cdots)$ in the case of the standard $m$-growth conditions was also proved working (e.g. \cite{Gianquintal3, Gianquintal2, Ladyzhenskaya, Morrey68}). This method is also available to investigate the variational problems and the regularity of solutions of quasi-linear elliptic equations. The regularity of solutions under nonstandard growth conditions is investigated  following an counterexample given by Giaquinta in 1987 \cite{Giaquinta1}, further relevant contributions are for example in \cite{AcerbiFusco94, AcerbiFusco94a, Marcellini89, Marcellini91, Marcellini93, Lieberman91, Lieberman93, Lieberman94, Mingione01, Lieberman05, Schmidt08, Schmidt09} and a more recent paper \cite{Hasto17}.

Our aim in the current paper is to study the H\"{o}lder continuity of the
minimizers for functionals defined on the Musielak-Orlicz-Sobolev
spaces and the H\"{o}lder continuity of the weak solutions for the associated
fully nonlinear elliptic equations, mainly in the framework of Musielak-Orlicz-Sobolev space and from the viewpoint of PDEs.

In fact, there are some important classical regularity results for the minimizers of integral functionals within the Sobolev, variable exponent, Orlicz and Musielak-Orlicz-Sobolev framework settings in the literature. 

In an early work of \cite{AcerbiFusco89}, Acerbi and Fusco proved that for any $W_{\text{loc}}^{1,p}(\Omega;\mathbb{R}^N)$-local minimizer $u$ with $1<p<2$ of the integral functional $\int f(x,v(x),Dv(x))\,\mathrm{d}x$,  its gradient $Du$ is actually locally $\lambda$-H\"{o}lder continuous for some $\lambda>0$ when $f$ fulfills some uniformly $p$-exponent increasing conditions. In the recent papers \cite{Mingione15, Mingione15s}, Colombo and  Mingione investigated the regularity of $W^{1,p}(\Omega;\mathbb{R}^N)$-local minimizer $u$ of the integral functional defined by $\mathcal{P}_{p,q}(w,\Omega)=\int_{\Omega}(|Dw|^p+a(x)|Dw|^q)\,\mathrm{d}x$ where $0\leq a(x)\leq L$ and $1<p<q$. In fact, the authors proved that $Du$ is actually locally H\"{o}lder continuous when $0\leq a(\cdot)\in C^{0,\alpha}(\Omega)$ and $\frac{q}{p}<1+\frac{\alpha}{n}$. The results in \cite{Mingione15} cover more general functionals than $\mathcal{P}_{p,q}$ in \cite{Mingione15s}. As described in \cite{Mingione15}, the result was actually extended to the vector case and a larger class of more general functionals. When considering the particular case $p=q$ or $a(x)\equiv0$, the regularity theory of minimizers is by now well understood, see for instance \cite{Manfredi88, Manfredi86, Mingione06}. In a much more recent paper \cite{Mingione16}, the authors considered the minimizer $u\in W^{1,1}(\Omega)$ of the functional defined by $\mathcal{P}_{\log}(w,\Omega):=\int_{\Omega}[|Dw|^p+a(x)|Dw|^p\log(e+|Dw|)]\,\mathrm{d}x$, in which the function $a(\cdot)$ is nonnegative, bounded and satisfies $|a(x)-a(y)|\leq w(|x-y|)$ for every $x,y \in\Omega$. They proved that if $l:=\limsup_{r\longrightarrow 0}w(r)\log(\frac{1}{r})<\infty$, then $u\in C^{1,\beta}_{loc}$ for some $\beta\in(0,1)$; if $l=0$, then $u\in C^{1,\beta}_{loc}$ for every $\beta\in(0,1)$; if $w(r)\lesssim r^{\sigma}$ with $\sigma\in(0,1)$ then $Du$ is locally H\"{o}lder continuous in $\Omega$. 

We point out that some regularity results in the variable exponent spaces framework can be found in the works \cite{Mingione99, Fan_Zhao, DieningHasto09}.  In the paper \cite{Mingione99}, Coscia and Mingione proved that for any $W^{1,1}(\Omega)$-local minimizer $u$ of the integral functional $\int_{\Omega}|\nabla u|^{p(x)}\,\mathrm{d}x$, its gradient $Du$ is actually locally H\"{o}lder continuous when $p(\cdot)$ is locally H\"{o}lder continuous in $\Omega$. When $p(\cdot)$ satisfies $R^{-\text{osc}\{p;B_R\}}\leq L,$ for all $B_{R}\subset\Omega$, Fan and Zhao \cite{Fan_Zhao} proved that $W^{1,p(\cdot)}$-minimizer $u\in C^{0,\beta}(\Omega)$ for some $\beta\in(0,1)$. In \cite{DieningHasto09}, Diening and H\"{a}st\"{o} also introduced in the important study \cite{DieningHasto09} the Triebel-Lizorkin spaces with variable smoothness and itegrability, including a trace theorem in the variable index case. In \cite{Adamowicz15}, Adamowicz and collaborators showed the continuity of quasiminimizers of energy functionals $\int f(x,u,\nabla u)\,\mathrm{d}x$ when $f$ satisfies some uniformly $p(\cdot)$-exponent monotonicity assumptions. 

In the recent works of \cite{Diening09, Diening11}, Diening and his collaborators proved a series of regularity results in Orlicz spaces. More precisely, they proved in \cite{Diening09} the $C^{1,\alpha}$-regularity for local minimizers of functionals with $\varphi$-growth including the decay estimate, where $\varphi$ is a convex $C^1$-function independent of the parameter $x\in\Omega\subset\mathbb{R}^N$; in \cite{Diening11}, they established a local Lipschitz result for the local minimizers of asymptotically convex variational integrals.

For regularity results in the Musielak-Orlicz-Sobolev setting, we also noticed that H\"{a}st\"{o} and collaborators proved, in an important recent work \cite{Hasto17}, that \textit{Harnack's inequality} still holds for quasi-minimizers in the Musielak-Orlicz-Sobolev spaces without any polynomial growth or coercivity conditions, which yields the local H\"{o}lder continuity of quasi-minimizers. Comparing with our current study, it is interesting to notice that we proposed a different monotonicity assumption for the $\Phi$ function from that in \cite{Hasto17}'s. In our current study, we have proposed a more general uniformly monotonicity condition on the $N(\Omega)$ function. Meanwhile, with the regularity results in the key Theorem \ref{key thm}, we can prove \textit{not only} the H\"{o}lder continuity of the minimizers for a more general class of energy functionals (see Section \ref{Sec4}), \textit{but also} the H\"{o}lder continuity of a kind of weak solutions for a class of fully nonlinear elliptic equations (see Section \ref{Sec5}).

The paper is organized as follows. In Section \ref{Sec2}, for the
readers' convenience we recall some definitions and properties about
Musielak-Orlicz-Sobolev spaces. In Section \ref{Sec3}, we give some crucial lemmas
in order to prove the main theorems of this paper. In Section \ref{Sec4},
we prove the H\"{o}lder continuity of the minimizers of a class of the energy functionals
in Musielak-Orlicz-Sobolev spaces. In Section \ref{Sec5}, the H\"{o}lder continuity of the
weak solution to a class of fully nonlinear elliptic equations is provided.

\section{The Musielak-Orlicz-Sobolev Spaces}\label{Sec2}

In this section, we list some definitions and propositions related to
Musielak-Orlicz-Sobolev spaces.
Firstly, we give the definition of \textit{$N$-function} and
\textit{generalized $N$-function} as following.

\vspace{0.3cm}

\begin{defi}
A function $A:\mathbb{R}\rightarrow[0,+\infty)$ is called an
$N$-function, denoted by $A\in N$, if $A$ is even and convex,
$A(0)=0, 0< A(t)\in C^0$ for $t\not=0$, and the following conditions
hold
\begin{equation*}
\lim_{t\rightarrow0+}\frac{A(t)}{t}=0\text{ and }
\lim_{t\rightarrow+\infty}\frac{A(t)}{t}=+\infty.
\end{equation*}
Let $\Omega$ be a smooth domain in $\mathbb{R}^n$. A function $A:\Omega\times\mathbb{R}\rightarrow[0,+\infty)$ is
called a generalized $ N$-function, denoted by $A\in N(\Omega)$, if
for each $t\in[0,+\infty)$, the function $A(\cdot,t)$ is measurable,
and for a.e. $x\in\Omega$, we have $A(x,\cdot)\in N$.
\end{defi}

\vspace{0.3cm}

Let $A\in N(\Omega)$, the Musielak-Orlicz space $L^{A}(\Omega)$ is
defined by
\begin{equation*}
\begin{aligned}
L^{A}(\Omega)&:=\bigg\{u:\,u\text{ is a measurable real function, and
}\exists\lambda>0\\
&\quad\quad\quad\quad\quad\quad\quad\quad\quad\quad\quad\quad\text{ such that
}\int_{\Omega}A\bigg(x,\frac{|u(x)|}{\lambda}\bigg)\,\mathrm{d}x<+\infty\bigg\}
\end{aligned}
\end{equation*}
with the (Luxemburg) norm
\begin{equation*}
\|u\|_{L^{A}(\Omega)}=\|u\|_A:=\inf\bigg\{\lambda>0:\,\int_{\Omega}A\bigg(x,\frac{|u(x)|}
{\lambda}\bigg)\,\mathrm{d}x\leq1\bigg\}.
\end{equation*}

The Musielak-Orlicz-Sobolev space $W^{1,A}(\Omega)$ can be defined by
\begin{equation*}
W^{1,A}(\Omega):=\{u\in L^{A}(\Omega):\,|\nabla u|\in
L^{A}(\Omega)\}
\end{equation*}
with the norm
\begin{equation*}
\|u\|_{W^{1,A}(\Omega)}=\|u\|_{1,A}:=\|u\|_A+\|\nabla u\|_{A},
\end{equation*}
where $\|\nabla u\|_{A}:=\|\,|\nabla u|\,\|_{A}$.

$A$ is called locally integrable if $A(\cdot,t_0)\in
 L_{\text{loc}}^1(\Omega)$ for every $t_0>0$.

\vspace{0.3cm}

\begin{defi} We say that $a(x,t)$ is the Musielak derivative of $A(x,t)\in N(\Omega)$ at $t$ if
for $x\in\Omega$ and $t\geq0$, $a(x,t)$ is the right-hand derivative of
$A(x,\cdot)$ at $t$; and for $x\in\Omega$ and $t\leq0$,
$a(x,t):=-a(x,-t)$.
\end{defi}

\vspace{0.3cm}

Define $\widetilde A:\Omega\times\mathbb{R}\rightarrow[0,+\infty)$
by
\begin{equation*}
\widetilde A(x,s)=\sup_{t\in\mathbb{R}}\big(st-A(x,t)\big)\text{ for
}x\in\Omega\text{ and }s\in\mathbb{R}.
\end{equation*}
$\widetilde A$ is called the \textit{complementary function} to $A$ in the
sense of Young. It is well known that if $A\in N(\Omega)$, then
$\widetilde A\in N(\Omega)$ and $A$ is also the complementary
function to $\widetilde A$.

For $x\in\Omega$ and $s\geq0$, we denote by $a_+^{-1}(x,s)$ the
right-hand derivative of $\widetilde{A}(x,\cdot)$ at $s$ at the same time
define $a_+^{-1}(x,s)= -a_+^{-1}(x,-s)$ for $x\in\Omega$ and
$s\leq0$. Then for $x\in\Omega$ and $s\geq0$, we have
\begin{equation*}
a_+^{-1}(x,s)=\sup\{t\geq0:\, a(x,t)\leq s\}=\inf\{t>0:\,a(x,t)>s\}.
\end{equation*}

\vspace{0.3cm}

\begin{prop}\label{Aa} (See \cite{Fan1, Musielak}) Let $A\in N(\Omega)$. Then the
following assertions hold.
\begin{enumerate}
\item $A(x,t)\leq a(x,t)t\leq A(x,2t)$ for $x\in\Omega$ and
$t\in\mathbb{R}$;
\item $A$ and $\widetilde A$ satisfy the Young inequality
\begin{equation*}
st\leq A(x,t)+\widetilde A(x,s) \text{ for }x\in\Omega \text{ and }
s,t\in\mathbb{R}
\end{equation*}
and the equality holds if $s=a(x,t)$ or $t=a_+^{-1}(x,s)$.
\end{enumerate}
\end{prop}

\vspace{0.3cm}

Let $A,B\in N(\Omega)$. We say that $A$ is weaker than $B$, denoted
by $A\preccurlyeq B$, if there exist positive constants $K_1,K_2$
and $h\in L^1(\Omega)\cap L^\infty(\Omega)$ such that
\begin{equation}\label{preccurlyeq}
A(x,t)\leq K_1B(x,K_2t)+h(x)\text{ for }x\in\Omega\text{ and
}t\in[0,+\infty).
\end{equation}

\vspace{0.3cm}

\begin{prop}(See \cite{Fan1, Musielak})
Let $A, B\in N(\Omega)$ and $A\preccurlyeq B$. Then $\widetilde
B\preccurlyeq \widetilde A$, $L^B(\Omega)\hookrightarrow
L^A(\Omega)$ and $L^{\widetilde A}(\Omega)\hookrightarrow
L^{\widetilde B}(\Omega)$.
\end{prop}

\vspace{0.3cm}

\begin{defi}\label{delta_2}
We say that a function $A:[0,+\infty)\rightarrow[0,+\infty)$
satisfies the $\Delta_2(\Omega)$ condition, denoted by $A\in
\Delta_2(\Omega)$, if there exist a positive constant $K>0$ and a
nonnegative function $h\in L^1(\Omega)$ such that
\begin{equation*}
A(x,2t)\leq KA(x,t)+h(x)\text{ for }x\in\Omega\text{ and
}t\in[0,+\infty).
\end{equation*}
\end{defi}

If $A(x,t)=A(t)$ is an $N$-function and $h(x)\equiv0$ in $\Omega$ in
Definition \ref{delta_2}, then $A\in\Delta_2(\Omega)$ if and only if
$A$ satisfies the well-known $\Delta_2$ condition defined in
\cite{Adams, Trudinger_1971}.

\vspace{0.3cm}

\begin{prop}\label{int_A} (See \cite{Fan1})
Let $A\in N(\Omega)$ satisfy $\Delta_2(\Omega)$. Then the following
assertions hold,
\begin{enumerate}
\item $L^A(\Omega)=\{u:\,u \text{ is a measurable function, and
}\int_{\Omega}A(x,|u(x)|)\,\mathrm{d}x<+\infty\}$;
\item $\int_\Omega A(x,|u|)\,\mathrm{d}x<1 \text{ (resp. } =1; >1) \Longleftrightarrow \|u\|_A<1 \text{ (resp. } =1;
>1)$, where $u\in L^A(\Omega)$;
\item $\int_\Omega A(x,|u_n|)\,\mathrm{d}x\rightarrow0 \text{ (resp. } 1; +\infty) \Longleftrightarrow \|u_n\|_A\rightarrow0 \text{ (resp. }1;
+\infty)$, where $\{u_n\}\subset L^A(\Omega)$;
\item $u_n\rightarrow u$ in $L^A(\Omega)\Longrightarrow\int_\Omega \big|A(x,|u_n|)\,\mathrm{d}x-
A(x,|u|)\big|\,\mathrm{d}x\rightarrow0$ as $n\rightarrow\infty$;
\item If $A'$ also satisfies $(\Delta_2)$, then
\begin{equation*}
\bigg|\int_{\Omega}u(x)v(x)\,\mathrm{d}x\bigg|\leq2\|u\|_A\|v\|_{\widetilde
A},\Forall u\in L^A(\Omega),v\in L^{\widetilde A}(\Omega);
\end{equation*}
\item $a(\cdot,|u(\cdot)|)\in L^{\widetilde A}(\Omega)$ for every $u\in
L^A(\Omega)$.
\end{enumerate}
\end{prop}

\vspace{0.3cm}

The following assumptions will be used.

\begin{enumerate}
\item[$(C_1)$] $\inf_{x\in\Omega}A(x,1)=c_1>0$;
\end{enumerate}

\vspace{0.3cm}

\begin{prop}(See \cite{Fan1})
If $A\in N(\Omega)$ satisfies $(C_1)$, then
$L^A(\Omega)\hookrightarrow L^1(\Omega)$ and
$W^{1,A}(\Omega)\hookrightarrow W^{1,1}(\Omega)$.
\end{prop}

\vspace{0.3cm}

Let $A\in  N(\Omega)$ be locally integrable. We will denote
\begin{equation*}
\begin{aligned}
W_{0}^{1,A}(\Omega):&=\overline{C_0^{\infty}(\Omega)}^{\|\,\cdot\,\|_{W^{1,A}(\Omega)}}\\
\mathcal{D}_0^{1,A}(\Omega):&=\overline{C_0^{\infty}(\Omega)}^{\|\nabla\,\cdot\,\|_{L^{A}(\Omega)}}.
\end{aligned}
\end{equation*}
In the case that $\|\nabla u\|_{A}$ is an equivalent norm in
$W_0^{1,A}(\Omega)$,
$W_0^{1,A}(\Omega)=\mathcal{D}_0^{1,A}(\Omega)$.

\vspace{0.3cm}

\begin{prop}(See \cite{Fan1})
Let $A\in N(\Omega)$ be locally integrable and satisfy $(C_1)$. Then
\begin{enumerate}
\item the spaces
$W^{1,A}(\Omega),W_{0}^{1,A}(\Omega)$ and
$\mathcal{D}_0^{1,A}(\Omega)$ are separable Banach spaces, and
\begin{equation*}
\begin{aligned}
W_{0}^{1,A}(\Omega)&\hookrightarrow W^{1,A}(\Omega)\hookrightarrow
W^{1,1}(\Omega)\\
\mathcal{D}_0^{1,A}(\Omega)&\hookrightarrow
\mathcal{D}_0^{1,1}(\Omega)=W_{0}^{1,1}(\Omega);
\end{aligned}
\end{equation*}
\item the spaces $W^{1,A}(\Omega),W_{0}^{1,A}(\Omega)$ and
$\mathcal{D}_0^{1,A}(\Omega)$ are reflexive provided $L^A(\Omega)$
is reflexive.
\end{enumerate}
\end{prop}

\vspace{0.3cm}

\begin{prop}\label{Poincare} (See \cite{Fan1})
Let $A,B\in  N(\Omega)$ and $A$ be locally integrable. If there is a
compact imbedding $W^{1,A}(\Omega)\hookrightarrow\hookrightarrow
L^{B}(\Omega)$ and $A\preccurlyeq B$, then there holds the following
Poincar\'{e} inequality
\begin{equation*}
\|u\|_A\leq c\|\nabla u\|_A, \Forall u\in W_0^{1,A}(\Omega),
\end{equation*}
which implies that $\|\nabla\cdot\|_A$ is an equivalent norm in
$W_0^{1,A}(\Omega)$ and
$W_0^{1,A}(\Omega)=\mathcal{D}_0^{1,A}(\Omega)$.
\end{prop}

\vspace{0.3cm}

The following assumptions will be used.

\begin{enumerate}
\item[$(P_1)$] $\Omega\subset\mathbb{R}^n(n\geq2)$ is a bounded domain with the cone property,
and $A\in N(\Omega)$;
\item[$(P_2)$]
$A:\overline{\Omega}\times\mathbb{R}\rightarrow[0,+\infty)$ is
continuous and $A(x,t)\in(0,+\infty)$ for $x\in\overline{\Omega}$
and $t\in(0,+\infty)$.
\end{enumerate}

\vspace{0.3cm}

Let $A$ satisfy $(P_1)$ and $(P_2)$. Denote by $A^{-1}(x,\cdot)$ the
inverse function of $A(x,\cdot)$.  We always assume that the
following condition holds.
\begin{enumerate}
\item[$(P_3)$] $A\in N(\Omega)$ and
\begin{equation}\label{0_1}
\int_0^1\frac{A^{-1}(x,t)}{t^{\frac{n+1}{n}}}\d t<+\infty,\Forall
x\in\overline\Omega.
\end{equation}
\end{enumerate}

\vspace{0.3cm}

Under assumptions $(P_1)$, $(P_2)$ and $(P_3)$, for each
$x\in\overline{\Omega}$, the function
$A(x,\cdot):[0,+\infty)\rightarrow[0,+\infty)$ is a strictly
increasing homeomorphism. Define a function $A_*^{-1}:
\overline{\Omega}\times[0,+\infty)\rightarrow[0,+\infty)$ by
\begin{equation}\label{inversA_*}
A_*^{-1}(x,s)=\int_0^s\frac{A^{-1}(x,\tau)}{\tau^{\frac{n+1}{n}}}\,\mathrm{d}\tau\text{
for }x\in\overline{\Omega} \text{ and }s\in[0,+\infty).
\end{equation}
Then under the assumption $(P_3)$, $A_*^{-1}$ is well defined, and
for each $x\in\overline{\Omega}$, $A_*^{-1}(x,\cdot)$ is strictly
increasing, $A_*^{-1}(x,\cdot)\in C^1((0,+\infty))$ and the function
$A_*^{-1}(x,\cdot)$ is concave.

Set
\begin{equation}\label{T}
T(x)=\lim_{s\rightarrow+\infty}A_*^{-1}(x,s), \Forall
x\in\overline\Omega.
\end{equation}
Then $0<T(x)\leq +\infty$. Define an even function $A_*:
\overline{\Omega}\times\mathbb{R}\rightarrow[0,+\infty)$ by
\begin{equation*}
\begin{aligned}
A_*(x,t)=\left\{ \begin{array}{ll}
          s,  & \text{ if } x\in \overline{\Omega}, |t|\in[0,T(x))\text{ and }A_*^{-1}(x,s)=|t|,\\
          +\infty,   & \text{ for } x\in \overline{\Omega} \text{ and } |t|\geq T(x).
                \end{array}\right.
\end{aligned}
\end{equation*}
Then if $A\in N(\Omega)$ and $T(x)=+\infty$ for any
$x\in\overline{\Omega}$, it is well known that $A_*\in N(\Omega)$
(see \cite{Adams}). $A_*$ is called the Sobolev conjugate function
of $A$ (see \cite{Adams} for the case of Orlicz functions).

Let $X$ be a metric space and $f:X\rightarrow(-\infty,+\infty]$ be
an extended real-valued function. For $x\in X$ with $f(x)\in
\mathbb{R}$, the continuity of $f$ at $x$ is well defined. For $x\in
X$ with $f(x)=+\infty$, we say that $f$ is continuous at $x$ if
given any $M>0$, there exists a neighborhood $U$ of $x$ such that
$f(y)>M$ for all $y\in U$. We say that
$f:X\rightarrow(-\infty,+\infty]$ is continuous on $X$ if $f$ is
continuous at every $x\in X$. Define Dom$(f)=\{x\in X :
f(x)\in\mathbb{R}\}$ and denote by $C^{1-0}(X)$ the set of all
locally Lipschitz continuous real-valued functions defined on $X$.

\vspace{0.3cm}

\begin{rem}\label{rem}
Suppose that $A\in N(\Omega)$ satisfy $(P_2)$, then for each $t_0\geq0, \widetilde{A}(x,t_0),\\A_*(x,t_0)$ are bounded.
\end{rem}

\vspace{0.3cm}

The following assumptions will also be used.
\begin{enumerate}
\item[$(P_4)$] $T:\overline{\Omega}\rightarrow[0,+\infty]$ is
continuous on $\overline{\Omega}$ and $T\in C^{1-0}(\text{Dom}(T))$;
\item[$(P_5)$]
$A_*\in C^{1-0}(\text{Dom}(A_*))$ and there exist three positive constants
$\delta_0$, $C_0$ and $t_0$ with $\delta_0<\frac{1}{n}$,
$0<t_0<\min_{x\in\overline{\Omega}}T(x)$ such that
\begin{equation*}
|\nabla_x A_*(x,t)|\leq
C_0(A_*(x,t))^{1+\delta_0},\quad{j=1,\dots,n},
\end{equation*}
for $x\in\Omega$ and $|t|\in[t_0,T(x))$ provided $\nabla_x A_*(x,t)$
exists.
\end{enumerate}

Let $A,B\in N(\Omega)$. We say that $A\ll B$ if, for any $k
> 0$,
\begin{equation*}
\lim_{t\rightarrow+\infty}\frac{A(x,kt)}{B(x,t)}=0\text{ uniformly
for }x\in\Omega.
\end{equation*}

\vspace{0.3cm}

\begin{rem}\label{Symbol}
Suppose that $A,B\in N(\Omega)$, then $A\ll B\Rightarrow
A\preccurlyeq B$.
\end{rem}

\vspace{0.3cm}

Next we give two embedding theorems for Musielak-Orlicz-Sobolev spaces recently
developed by Fan in \cite{Fan2}.

\vspace{0.3cm}

\begin{thm}\label{imbedding}(See \cite{Fan2}, \cite{Liu_Zhao_15})
Let $(P_1)-(P_5)$ hold. Then
\begin{enumerate}
\item[(i)] There is a continuous imbedding $W^{1,A}(\Omega)\hookrightarrow
L^{A_*}(\Omega)$;
\item[(ii)] Suppose that $B\in N(\Omega)$,
$B:\overline{\Omega}\times[0,+\infty)\rightarrow[0,+\infty)$ is
continuous, and $B(x,t)\in(0,+\infty)$ for $x\in\Omega$ and
$t\in(0,+\infty)$. If $B\ll A_*$, then there is a compact imbedding
$W^{1,A}(\Omega)\hookrightarrow\hookrightarrow L^B(\Omega)$.
\end{enumerate}
\end{thm}

\vspace{0.3cm}

By Theorem \ref{imbedding}, Remark \ref{Symbol} and Proposition
\ref{Poincare}, we have the following:

\vspace{0.3cm}

\begin{thm}\label{embeddings}(See \cite{Fan2}, \cite{Liu_Zhao_15})
Let $(P_1)-(P_5)$ hold and furthermore, $A,A_*\in N(\Omega)$. Then
\begin{enumerate}
\item[(i)] $A\ll A_*$, and there is a compact imbedding $W^{1,A}(\Omega)\hookrightarrow\hookrightarrow
L^A(\Omega)$;
\item[(ii)] there holds the poincar\'{e}-type inequality
\begin{equation*}
\|u\|_A\leq C\|\nabla u\|_A\text{ for }u\in W_0^{1,A}(\Omega),
\end{equation*}
i.e. $\|\nabla u\|_A$ is an equivalent norm on $W_0^{1,A}(\Omega)$.
\end{enumerate}
\end{thm}

\vspace{0.3cm}

\section{Some Lemmas}\label{Sec3}
Suppose $\Omega\subset\mathbb{R}^n$ is a bounded smooth domain, and $A\in N(\Omega)$
satisfies the following Condition $(\mathscr{A})$, denoted by
$A\in\mathscr{A}$.
\begin{enumerate}
\item[$(\mathscr{A})$]
$A\in N(\Omega)$ satisfies assumptions $(P_1)$, $(P_2)$, $(P_3)$, $(P_5)$ in Section
\ref{Sec2} and the following
\begin{enumerate}
\item[$(\widetilde{P_4})$] $T(x)$ defined in \eqref{T} satisfies
$T(x)=+\infty$ for all $x\in\overline{\Omega}$.
\end{enumerate}
\end{enumerate}

\vspace{0.3cm}

\begin{lem}\label{compare}(see \cite{liuandyao})
Suppose that $A\in N(\Omega)$, and there exists a strictly
increasing differentiable function
$\mathfrak{A}:[0,+\infty)\rightarrow[0,+\infty)$ such that
\begin{equation}\label{A_A}
A(x,\alpha t)\geq \mathfrak{A}(\alpha)A(x,t), \,\forall\,
\alpha\geq0,t\in\mathbb{R},x\in\Omega.
\end{equation}
\begin{enumerate}
\item[(i)] Then there exists a
strictly increasing differentiable function
$\widehat{\mathfrak{A}}:[0,+\infty)\rightarrow[0,+\infty)$, defined
by
\begin{equation}\label{special}
\begin{aligned}
\widehat{\mathfrak{A}}(\beta)=\left\{ \begin{array}{ll}
          \frac{1}{\mathfrak{A}(\frac{1}{\beta})},  & \text{ for } \beta>0, \\
          0,   & \text{ for } \beta=0,
                \end{array}\right.
\end{aligned}
\end{equation}
such that
\begin{equation}\label{A_star_1}
A(x,\beta t)\leq \widehat{\mathfrak{A}}(\beta)A(x,t), \,\forall
\beta>0,t\in\mathbb{R},x\in\Omega,
\end{equation}
and furthermore $\widehat{\widehat{\mathfrak{A}}}=\mathfrak{A}$;
\item[(ii)] If $\mathfrak{A}$ satisfies
\begin{equation}\label{plessN}
n\mathfrak{A}(\alpha)>\alpha\mathfrak{A}'(\alpha),
\end{equation}
then $A_*\in N(\Omega)$, and there exists a
strictly increasing differentiable function
$\mathfrak{A}_*:[0,+\infty)\rightarrow[0,+\infty)$, defined by
\begin{equation}
\begin{aligned}\label{exprr2}
\mathfrak{A}_*^{-1}(\sigma)=\left\{ \begin{array}{ll}
          \frac{1}{\sigma^{\frac{1}{n}}\mathfrak{A}^{-1}(\sigma^{-1})},  & \text{ for } \sigma>0, \\
          0,   & \text{ for } \sigma=0,
                \end{array}\right.
\end{aligned}
\end{equation}
such that
\begin{equation}\label{A_star}
A_*(x,\beta t)\leq \mathfrak{A}_*(\beta)A_*(x,t), \,\forall
\beta>0,t\in\mathbb{R},x\in\Omega;
\end{equation}
\item[(iii)] If $\mathfrak{A}$ satisfies
\begin{equation}\label{plessN2}
\alpha\mathfrak{A}'(\alpha)>\mathfrak{A}(\alpha),
\end{equation}
then $\widetilde A\in N(\Omega)$, and there exists two
strictly increasing differentiable functions
$\widetilde{\mathfrak{A}},\widehat{\widetilde{\mathfrak{A}}}:[0,+\infty)\rightarrow[0,+\infty)$, defined by
\begin{equation}\label{exprr3}
\begin{aligned}
\widetilde{\mathfrak{A}}^{-1}(\sigma)=\left\{ \begin{array}{ll}
         \frac {\sigma}{\mathfrak{A}^{-1}(\sigma)},  & \text{ for } \sigma>0, \\
          0,   & \text{ for } \sigma=0,
                \end{array}\right.
\end{aligned}
\end{equation}
and
\begin{equation*}
\begin{aligned}
\widehat{\widetilde{\mathfrak{A}}} ^{-1}(\sigma)=\left\{ \begin{array}{ll}
         \sigma\mathfrak{A}^{-1}(\sigma^{-1}),  & \text{ for } \sigma>0, \\
          0,   & \text{ for } \sigma=0,
                \end{array}\right.
\end{aligned}
\end{equation*}
such that
\begin{equation}\label{A_star2}
\widetilde A(x,\beta t)\leq \widetilde{\mathfrak{A}}(\beta)\widetilde A(x,t), \,\forall
\beta>0,t\in\mathbb{R},x\in\Omega.
\end{equation}
\begin{equation}\label{A_star3}
\widetilde A(x,\beta t)\geq \widehat{\widetilde{\mathfrak{A}}} (\beta)\widetilde A(x,t), \,\forall
\beta>0,t\in\mathbb{R},x\in\Omega.
\end{equation}
\end{enumerate}
\end{lem}

\vspace{0.3cm}

\begin{defi}
We say that $\mathfrak{C}:\mathbb{R}^+\rightarrow\mathbb{R}^+$
satisfies Condition $\Delta_{\mathbb{R}^+}$, denoted by
$\mathfrak{C}\in\Delta_{\mathbb{R}^+}$, if there exists a constant
$M_0>0$ such that
\begin{equation}\label{multi}
\mathfrak{C}(\alpha\beta)\leq
M_0\mathfrak{C}(\alpha)\mathfrak{C}(\beta),\,\forall\,\alpha,\beta>0.
\end{equation}
\end{defi}

\vspace{0.3cm}

\begin{rem}\label{rem1}
By \eqref{special}, \eqref{exprr2} and \eqref{exprr3}, if $\,\mathfrak{A},\mathfrak{A}^{-1}\in\Delta_{\mathbb{R}^+}$, then $\widehat{\mathfrak{A}},\mathfrak{A}_*^{-1},\widetilde{\mathfrak{A}}^{-1},\widehat{\widetilde{\mathfrak{A}}}^{-1}\in\Delta_{\mathbb{R}^+}$, and there exists a constant $M_0>0$ such that the following two inequalities hold
\begin{equation*}
\mathfrak{A}^{-1}(\alpha)\mathfrak{A}^{-1}(\alpha^{-1})\leq M_0,\quad
\alpha\mathfrak{A}^{-1}(\widehat{\mathfrak{A}}(\alpha^{-1}))\leq M_0,\, \forall \alpha>0.
\end{equation*}
\end{rem}

\vspace{0.3cm}

\begin{lem}\label{sequence}(see \cite{liuandyao})
Let $\{y_h\}\subset\mathbb{R}$ be a sequence satisfying
\begin{equation}\label{n_n-1}
y_{h+1}\leq\frac{1}{\beta}\mathfrak{A}_*\bigg(\frac{\mathfrak{A}^{-1}(\beta)}{\beta^{\frac{1}{n}}}
c2^{h}\mathfrak{A}^{-1}(c\mathfrak{A}_*(2^{h+2})y_h)\bigg),\,\forall\,\beta>0,
\end{equation}
where $c$ is a positive constant. If
$\,\mathfrak{A},\mathfrak{A}^{-1},\mathfrak{A}_*\in\Delta_{\mathbb{R}^+}$
then there exists a $y_0^*>0$ such that for $y_0\leq y_0^*$, $y_h\rightarrow0$ as
$h\rightarrow\infty$.
\end{lem}

\vspace{0.3cm}

For a measurable set $E\subset\mathbb{R}^n$, we denote by $\text{mes}(E)$ or $|E|$ the $n$-Lebesgue measure of $E$. For a measurable function $u$ defined in $\Omega$ and a measurable set $E\subset\Omega$ denote
\begin{equation*}
\begin{aligned}
&\max_{E}u(x):=\text{ess sup}_{x\in E}u(x),\quad \min_{E}u(x):=\text{ess inf}_{x\in E}u(x),\\
&\text{osc}(u;E):=\max_Eu(x)-\min_Eu(x).
\end{aligned}
\end{equation*}
If $u\in W^{1,A}(\Omega)$ and
$B_\rho=B_\rho(x):=\{y\in\mathbb{R}^n:\,|y-x|<\rho\}\subset\Omega$ is any given ball, we denote $\Omega_{k,\rho}:=\{x\in
B_\rho:u(x)>k\}$, where $k$ is a real number.

\vspace{0.3cm}

The following two lemmas will be used.

\vspace{0.3cm}

\begin{lem}\label{lemma3}(see Lemma 3.9 of Chapter 2 of \cite{Ladyzhenskaya})
For any $u\in W^{1,1}(B_{\rho})$ and arbitrary number $k$ and $l$ with $l>k$, the following inequality holds
\begin{equation}\label{lem3.1}
(l-k)|\Omega_{l,\rho}|^{1-\frac{1}{n}}\leq C\frac{\rho^{n}}{|B_{\rho}-\Omega_{k,\rho}|}\int_{\Omega_{k,\rho}\setminus\Omega_{l,\rho}}|\nabla u|\,\mathrm{d}x,
\end{equation}
where $C=C(n)>1$ is a constant depending only on $n$.
\end{lem}

\vspace{0.3cm}

\begin{lem}\label{leylemma}(see Lemma 4.8 of Chapter 2 of \cite{Ladyzhenskaya})
Suppose a function $u(x)$ is measurable and bounded in some ball $B_{R_{0}}$. Consider balls $B_{R}$ and $B_{bR}$ which have a common center with $B_{R_{0}}$, where $b>1$ is a fixed constant. Suppose in addition that for any $0<R\leq b^{-1}R_{0}$ at least one of the following two inequalities is valid
\begin{equation*}
\text{osc}\{u;B_{R}\}\leq c_{1}R^{\varepsilon},\quad\text{osc}\{u,B_{R}\}\leq\theta\text{osc}\{u;B_{bR}\},
\end{equation*}
 where $c_{1}>0$, $\varepsilon\leq1$ and $\theta<1$ are positive constants. Then $u\in C^{0,\alpha}(B_{R_{0}})$, where $\alpha=\min\{\varepsilon,-\log_{b}\theta\}$.
\end{lem}

\vspace{0.3cm}

In the following lemma $A\in N(\Omega)\cap\mathscr{A}$ satisfies
the following assumptions:
\begin{enumerate}
\item[$(A_1)$]  There exists a strictly increasing
differentiable function
$\mathfrak{A}:[0,+\infty)\rightarrow[0,+\infty)$ satisfying
\begin{equation}\label{plessN1}
n\mathfrak{A}(\alpha)>\alpha\mathfrak{A}'(\alpha)>\mathfrak{A}(\alpha)
\end{equation}
such that
\begin{enumerate}
\item[($A_{11}$)] $A(x,\alpha t)\geq
\mathfrak{A}(\alpha)A(x,t), \,\forall\,
\alpha\geq0,t\in\mathbb{R},x\in\Omega;$
\item[($A_{12}$)] $\mathfrak{A},\mathfrak{A}^{-1},\mathfrak{A}_*,\widehat{\widetilde{\mathfrak{A}}},\in\Delta_{\mathbb{R}^
+}$;
\end{enumerate}
\end{enumerate}

\vspace{0.3cm}

\begin{lem}\label{keylemma}(see \cite{liuandyao})
Let $A\in N(\Omega)\cap\mathscr{A}$ satisfy $(A_1), B\in N(\Omega)$ satisfy
$B\preccurlyeq A_*,$ and $u\in W^{1,A}(\Omega)$ satisfy for
any $B_R\subset\Omega,R\leq R_0$ and for all $\sigma\in(0,1)$ and any
$k\geq t_0$,
\begin{equation}\label{key}
\int_{\Omega_{k,\sigma R}}A(x,|\nabla u|)\,\mathrm{d}x\leq
\gamma\int_{\Omega_{k,R}}A_*\bigg(x,\frac{u-k}{(1-\sigma)R}\bigg)\,\mathrm{d}x
+\gamma_{1}\int_{\Omega_{k,R}}B(x,k)\,\mathrm{d}x.
\end{equation}
Then $u$ is locally bounded above in $\Omega$.
\end{lem}

\vspace{0.3cm}

\begin{defi}\label{keyde}
Let $M$, $\gamma$, $\gamma_{1}$ and $\delta$ are positive constants with $\delta\leq2$. We will say that a function $u$ belongs to class $\mathscr{B}(\Omega,M,\gamma,\gamma_{1},\delta)$ if $u\in W^{1,A}(\Omega)$, $\max_\Omega|u(x)|\leq M$ and the functions $w(x)=\pm u(x)$ satisfy the inequality
\begin{equation}\label{keyeq}
\int_{\Omega_{k,\sigma \rho}}A(x,|\nabla w|)\,\mathrm{d}x\leq
\gamma\int_{\Omega_{k,\rho}}A\bigg(x,\frac{w(x)-k}{(1-\sigma)\rho}\bigg)\,\mathrm{d}x
+\gamma_{1}|\Omega_{k,\rho}|
\end{equation}
for arbitrary $B_{\rho}\subset\Omega,\sigma\in(0,1)$ and such $k$
\begin{equation}\label{keycon}
k\geq\max_{B_{\rho}}w(x)-\delta M.
\end{equation}
\end{defi}

\vspace{0.3cm}

The main result of this section is the following:

\vspace{0.3cm}

\begin{thm}\label{key thm}
Let $A\in N(\Omega)\cap\mathscr{A}$ satisfy $(A_1)$ and $(C_1)$. Then $\mathscr{B}(\Omega,M,\gamma,\gamma_{1},\delta)\subset C^{0,\alpha}(\Omega)$, where the constant $\alpha\in(0,1)$ depends only on the parameters $n$, $A$, $\gamma$ and $\delta$, but it is independent of $\gamma_{1}$ and $M$.
\end{thm}

\vspace{0.3cm}

We are now on the position to prove Theorem \ref{key thm}. Let $u\in\mathscr{B}(\Omega,M,\gamma,\gamma_{1},\delta)$ in which $A$ satisfies the condition $(A_1)$ and $(C_1)$. Without loss of generality we may assume that $M\geq1$ and $\gamma\geq1$.
To prove Theorem \ref{key thm} it suffices to prove that for each $x_{0}\in\Omega$ there is a ball $B_{R_{0}}\subset\Omega$ such that $u\in C^{0,\alpha}(B_{R_{0}})$ where $\alpha=\alpha(n,A,\gamma,\delta)$ is a constant.

Now let $x_{0}\in\Omega$ be given arbitrarily. Choose a positive number $R_{0}<1$ such that $\overline{B_{R_{0}}(x_{0})}\subset\Omega$, take arbitrarily $R\in(0,R_{0}]$. It is easy to see that at least one of the two functions $w=\pm u$ satisfies the following condition
\begin{equation}\label{cod1}
\text{mes}\bigg\{x\in B_{\frac{R}{2}}:w(x)>\max_{B_{R}}w(x)-\frac{1}{2}\text{osc}\{u;B_{R}\}\bigg\}\leq\frac{1}{2}\text{mes} B_{\frac{R}{2}}.
\end{equation}

From now on we denote by $w$ the function identified to $u$ or $-u$ that satisfies \eqref{cod1}. Set
\begin{equation}\label{cod2}
\tau=\max\{2,2/\delta\},\quad \psi=\tau^{-1}\text{osc}\{u;B_{R}\},\quad k'=\max_{B_{R}}w(x)-\psi.
\end{equation}
Then
\begin{equation}\label{cod3}
k'\geq\max_{B_{R}}w(x)-\frac{1}{2}\text{osc}\{u;B_{R}\}
\end{equation}
and
\begin{equation}\label{cod4}
k'\geq\max_{B_{R}}w(x)-\delta M.
\end{equation}

Note that \eqref{cod4} implies that \eqref{keyeq} and \eqref{keycon} hold for $k\geq k'$ and $\rho\leq R$.

Under the above assumptions we can give the proof of Theorem \ref{key thm} throughout the following Lemmas \ref{keylemma1}-\ref{keylemma4}.

\vspace{0.3cm}

\begin{lem}\label{keylemma1}
Let $A\in N(\Omega)\cap\mathscr{A}$ satisfy $(A_1)$ and $(C_1)$. Then there is a positive constant $\theta=\theta(n,A,\gamma)$ such that the following equation
\begin{equation}\label{key1.1}
|\Omega_{k^{0},\frac{R}{2}}|\leq\theta R^{n},
\end{equation}
implies that at least one of the following two inequalities holds:
\begin{equation}\label{key1.2}
H\leq\widehat{\mathfrak{A}}^{-1}\bigg(\frac{\gamma_{1}}{\gamma}\bigg) R,
\end{equation}
\begin{equation}\label{key1.3}
\max_{B_{R/4}}w(x)\leq k^{0}+\frac{H}{2},
\end{equation}
where
\begin{equation}\label{key1.4}
0<H<\psi, k^{0}=\max_{B_{R}}w(x)-H.
\end{equation}
\end{lem}

\vspace{0.3cm}

\begin{proof}

Set
\begin{equation*}
t=A^{-1}\bigg(x,\frac{s}{\mathfrak{A}(\alpha)}\bigg)
\end{equation*}
in \eqref{A_A}. Then for $\alpha>0$ we have
\begin{equation*}
A\bigg(x,\alpha
A^{-1}\bigg(x,\frac{s}{\mathfrak{A}(\alpha)}\bigg)\bigg)\geq
\mathfrak{A}(\alpha)A\bigg(x,A^{-1}\bigg(x,\frac{s}{\mathfrak{A}(\alpha)}\bigg)\bigg)=s,
\end{equation*}
or equivalently
\begin{equation*}
\alpha A^{-1}\bigg(x,\frac{s}{\mathfrak{A}(\alpha)}\bigg)\geq
A^{-1}(x,s),
\end{equation*}
which implies that
\begin{equation*}
\frac{\alpha}{(\mathfrak{A}(\alpha))^{\frac{1}{n}}}
\frac{A^{-1}\big(x,\frac{s}{\mathfrak{A}(\alpha)}\big)}
{\big(\frac{s}{\mathfrak{A}(\alpha)}\big)^{\frac{n+1}{n}}}\frac{1}{\mathfrak{A}(\alpha)}
\geq\frac{A^{-1}(x,s)}{s^{\frac{n+1}{n}}}.
\end{equation*}

Integrating the above inequality with respect to $s$ from $0$ to $t$, we have
\begin{equation*}
\frac{\alpha}{(\mathfrak{A}(\alpha))^{\frac{1}{n}}}\int_0^{\frac{t}{\mathfrak{A}(\alpha)}}
\frac{A^{-1}(x,r)}{r^{\frac{n+1}{n}}}\,\mathrm{d}r\geq\int_0^t
\frac{A^{-1}(x,s)}{s^{\frac{n+1}{n}}}\,\mathrm{d}s,
\end{equation*}
where $r=\frac{s}{\mathfrak{A}(\alpha)}$. Then the definition of
$A_*$ in \eqref{inversA_*} yields that
\begin{equation}\label{proccess}
\frac{\alpha}{(\mathfrak{A}(\alpha))^{\frac{1}{n}}}A_*^{-1}\bigg(x,\frac{t}{\mathfrak{A}(\alpha)}\bigg)
\geq A_*^{-1}(x,t).
\end{equation}

Setting $t=\mathfrak{A}(\alpha)A_*\bigg(x,\frac{\mathfrak{A}(\alpha)}{\alpha}\bigg)$ in \eqref{proccess}, we conclude that
\begin{equation*}
(\mathfrak{A}(\alpha))^{1-\frac{1}{n}}\geq A_*^{-1}\bigg(x,\mathfrak{A}(\alpha)A_*\bigg(x,\frac{\mathfrak{A}(\alpha)}{\alpha}\bigg)\bigg),
\end{equation*}
or equivalently
\begin{equation}\label{proccess1}
A_*(x,\mu^{1-\frac{1}{n}})\geq\mu A_*\bigg(x,\frac{\mu}{\mathfrak{A}^{-1}(\mu)}\bigg),
\end{equation}
in which $\mu=\mathfrak{A}(\alpha)$.
Then by \eqref{A_star} and the above \eqref{proccess1}, we have
\begin{equation*}
A_*(x,\mu^{1-\frac{1}{n}})\geq\mu A_*\bigg(x,\frac{\mu}{\mathfrak{A}^{-1}(\mu)}\bigg)\geq A_*\bigg(x,\mathfrak{A}_*^{-1}(\mu)\frac{\mu}{\mathfrak{A}^{-1}(\mu)}\bigg),
\end{equation*}
which implies that
\begin{equation}\label{Imeq}
\mathfrak{A}_*^{-1}(\mu)\mu\leq\mu^{1-\frac{1}{n}}\mathfrak{A}^{-1}(\mu).
\end{equation}

For $h=0,1,2\cdots $, set
\begin{equation*}
\begin{aligned}
 &\rho_{h}=\frac{R}{4}+\frac{R}{2^{h+2}}, &&k_{h}=k^{0}+\frac{H}{2}-\frac{H}{2^{h+1}},\\
 &y_{h}=R^{-n}|\Omega_{k_{h},\rho_{h}}|, &&D_{h+1}=\Omega_{k_{h},\rho_{h+1}}\setminus \Omega_{k_{h+1},\rho_{h+1}}.
\end{aligned}
\end{equation*}

Applying inequality \eqref{lem3.1} to $l=k_{h+1},k=k_{h}$ and $\rho=\rho_{h+1}$ we get
\begin{equation}\label{key1.5}
\begin{aligned}
&y_{h+1}^{1-\frac{1}{n}}\\
=&R^{1-n}|\Omega_{k_{h+1},\rho_{h+1}}|^{1-\frac{1}{n}}\\
\leq&\frac{CR^{1-n}}{(k_{h+1}-k_{h})}\frac{\rho^{n}_{h+1}}{|B_{\rho_{h+1}}-\Omega_{k_{h},\rho_{h+1}}|}
\int_{D_{h+1}}|\nabla w|\,\mathrm{d}x\\
\leq& CR^{1-n}2^{h+2}H^{-1}\bigg(\frac{R}{2}\bigg)^{n}\frac{1}{|B_{R/4}-\Omega_{k^{0},\frac{R}{2}}|}2|\nabla w|_{A;D_{h+1}}\cdot|1|_{\widetilde{A};D_{h+1}}\\
=&\frac{C2^{h+3-n}H^{-1}R}{|B_{R/4}-\Omega_{k^{0},\frac{R}{2}}|}\mathfrak{A}^{-1}\bigg(\mathfrak{A}(|\nabla w|_{A;D_{h+1}})\int_{D_{h+1}}A\bigg(x,\frac{|\nabla w|}{|\nabla w|_{A;D_{h+1}}}\bigg)\,\mathrm{d}x\bigg)\\
&\quad\quad\cdot\widehat{\widetilde{\mathfrak{A}}}^{-1}\bigg(\widehat{\widetilde{\mathfrak{A}}}(|1|_{\widetilde{A};D_{h+1}})\cdot
\int_{D_{h+1}}\widetilde{A}\bigg(x,\frac{1}{|1|_{\widetilde{A};D_{h+1}}}\bigg)\,\mathrm{d}x\bigg)\\
\leq& \frac{C2^{h+3-n}H^{-1}R}{|B_{R/4}-\Omega_{k^{0},\frac{R}{2}}|}\mathfrak{A}^{-1}\bigg(\int_{D_{h+1}}A(x,|\nabla w|)\,\mathrm{d}x\bigg)\cdot\widehat{\widetilde{\mathfrak{A}}}^{-1}\bigg(\int_{D_{h+1}}\widetilde{A}(x,1)\,\mathrm{d}x\bigg)\\
\leq& \frac{C2^{h+3-n}H^{-1}R}{|B_{R/4}-\Omega_{k^{0},\frac{R}{2}}|}\mathfrak{A}^{-1}\bigg(\int_{D_{h+1}}A(x,|\nabla w|)\,\mathrm{d}x\bigg)\cdot\widehat{\widetilde{\mathfrak{A}}}^{-1}(c_1|D_{h+1}|)\\
\leq&\frac{C2^{h+3-n}H^{-1}R}{|B_{R/4}-\Omega_{k^{0},\frac{R}{2}}|}\mathfrak{A}^{-1}\bigg(\int_{D_{h+1}}A(x,|\nabla w|)\,\mathrm{d}x\bigg)\cdot\widehat{\widetilde{\mathfrak{A}}}^{-1}(c_2R^ny_{h}).
\end{aligned}
\end{equation}

If \eqref{key1.1} holds, choose
\begin{equation*}
\theta\leq\frac{1}{2}\cdot4^{-n}\omega_{n}.
\end{equation*}
Then
\begin{equation*}
|B_{R/4}-\Omega_{k^{0},\frac{R}{2}}|\geq 4^{-n}R^{n}\omega_{n}-\frac{1}{2}\cdot4^{-n}\omega_{n}R^{n}=2^{-2n-1}\omega_{n}R^{n},
\end{equation*}
and consequently from \eqref{key1.5} it follows that
\begin{equation}\label{key1.7}
y_{h+1}^{1-\frac{1}{n}}\leq c2^{h}R^{1-n}H^{-1}\omega_{n}^{-1}\mathfrak{A}^{-1}\bigg(\int_{D_{h+1}}A(x,|\nabla w|)\,\mathrm{d}x\bigg)\widehat{\widetilde{\mathfrak{A}}}^{-1}(c_2R^ny_{h}).
\end{equation}
Applying inequality \eqref{keyeq} with $k=k_{h} , \rho=\rho_{h},\sigma \rho=\rho_{h+1}$, we obtain
\begin{equation}\label{key2}
\begin{aligned}
&\int_{D_{h+1}}A(x,|\nabla w|)\,\mathrm{d}x\\
\leq& \int_{\Omega_{k_{h},\rho_{h+1}}}A(x,|\nabla w|)\,\mathrm{d}x\\
\leq& \gamma\int_{\Omega_{k_{h},\rho_{h}}}A\bigg(x,\frac{w-k_{h}}{\rho_{h}-\rho_{h+1}}\bigg)\,\mathrm{d}x+\gamma_{1}|\Omega_{k_{h},\rho_{h}}|\\
\leq& \gamma\int_{\Omega_{k_{h},\rho_{h}}}A\bigg(x,\frac{2^{h+3}}{R}|w-k_{h}|\bigg)\,\mathrm{d}x+\gamma_{1}|\Omega_{k_{h},\rho_{h}}|\\
\leq& \gamma\int_{\Omega_{k_{h},\rho_{h}}}A\bigg(x,\frac{2^{h+3}}{R}H\bigg)\,\mathrm{d}x+\gamma_{1}|\Omega_{k_{h},\rho_{h}}|\\
\leq& \gamma\mathfrak{\widehat{A}}(R^{-1}H)\int_{\Omega_{k_{h},\rho_{h}}}A(x,2^{h+3})\,\mathrm{d}x+\gamma_{1}|\Omega_{k_{h},\rho_{h}}|\\
\leq& \gamma\mathfrak{\widehat{A}}(R^{-1}H)\bigg[\int_{\Omega_{k_{h},\rho_{h}}}A_{*}(x,2^{h+3})\,\mathrm{d}x+\int_{\Omega_{k_{h},\rho_{h}}}C_1\,\mathrm{d}x\bigg]
+\gamma_{1}|\Omega_{k_{h},\rho_{h}}|\\
\leq& \gamma\mathfrak{\widehat{A}}(R^{-1}H)[C_{2}\mathfrak{A}_{*}(2^{h+2})|\Omega_{k_{h},\rho_{h}}|+C_{1}|\Omega_{k_{h},\rho_{h}}|]
+\gamma_{1}|\Omega_{k_{h},\rho_{h}}|\\
\leq& C\gamma\mathfrak{\widehat{A}}(R^{-1}H)\mathfrak{A}_{*}(2^{h+2})|\Omega_{k_{h},\rho_{h}}|+\gamma_{1}|\Omega_{k_{h},\rho_{h}}|.
\end{aligned}
\end{equation}

Assume that equation \eqref{key1.2} does not hold, then we have
\begin{equation}\label{key3}
\int_{D_{h+1}}A(x,|\nabla w|)\,\mathrm{d}x\leq c\gamma\mathfrak{\widehat{A}}(R^{-1}H)\mathfrak{A}_{*}(2^{h+2})R^ny_h.
\end{equation}
Then from \eqref{key1.7}, we can get
\begin{equation}\label{key11.7}
y_{h+1}^{1-\frac{1}{n}}\leq c2^{h}R^{1-n}H^{-1}\omega_{n}^{-1}\mathfrak{A}^{-1}(c\gamma\mathfrak{\widehat{A}}(R^{-1}H)\mathfrak{A}_{*}(2^{h+2})R^ny_h)\widehat{\widetilde{\mathfrak{A}}}^{-1}(c_2R^ny_{h}).
\end{equation}

For any $\beta>0$, set $\mu=\beta y_{h+1}$ in \eqref{Imeq}. Then
\begin{equation*}
\mathfrak{A}_*^{-1}(\beta y_{h+1})\beta y_{h+1}\leq\beta^{1-\frac{1}{n}} y_{h+1}^{1-\frac{1}{n}}\mathfrak{A}^{-1}(\beta y_{h+1}),
\end{equation*}
or equivalently
\begin{equation}\label{Imequat}
\mathfrak{A}_*^{-1}(\beta y_{h+1}) y_{h+1}\leq\frac{1}{\beta^{\frac{1}{n}}} y_{h+1}^{1-\frac{1}{n}}\mathfrak{A}^{-1}(\beta y_{h+1}).
\end{equation}

By the definition of $y_h$, it is clear that $y_h\downarrow\varepsilon_0\geq0$ as $h\rightarrow\infty$. \\
\textbf{Claim: if \eqref{key1.2} does not hold, then $\varepsilon_0=0$.}

To prove $\varepsilon_0=0$, we argue by the contrary $\varepsilon_0>0$. Then Lemma \ref{compare}, Remark \ref{rem1}, $\mathfrak{A}^{-1},\mathfrak{A}\in\Delta_{\mathbb{R}^+}$, \eqref{key11.7} and \eqref{Imequat} imply
\begin{equation*}
\begin{aligned}
&\mathfrak{A}_*^{-1}(\beta y_{h+1}) y_{h+1}\\
\leq&\frac{1}{\beta^{\frac{1}{n}}} y_{h+1}^{1-\frac{1}{n}}\mathfrak{A}^{-1}(\beta y_{h+1})\\
\leq& \frac{\mathfrak{A}^{-1}(\beta)}{\beta^{\frac{1}{n}}} cy_{h+1}^{1-\frac{1}{n}}\mathfrak{A}^{-1}(y_{h+1})\\
\leq&\frac{\mathfrak{A}^{-1}(\beta)}{\beta^{\frac{1}{n}}}c2^{h}R^{1-n}H^{-1}\omega_{n}^{-1}\mathfrak{A}^{-1}(c\gamma\mathfrak{\widehat{A}}(R^{-1}H)\mathfrak{A}_{*}(2^{h+2})R^ny_h)\\
&\quad\quad\cdot\widehat{\widetilde{\mathfrak{A}}}^{-1}(c_2R^ny_{h})\mathfrak{A}^{-1}(y_{h+1})\\
=&\frac{\mathfrak{A}^{-1}(\beta)}{\beta^{\frac{1}{n}}}c2^{h}R^{1-n}H^{-1}\omega_{n}^{-1}\mathfrak{A}^{-1}(c\gamma\mathfrak{\widehat{A}}(R^{-1}H)\mathfrak{A}_{*}(2^{h+2})R^ny_h)\\
&\quad\quad\cdot c_2R^ny_{h}\mathfrak{A}^{-1}\bigg(\frac{1}{c_2R^ny_{h}}\bigg)\mathfrak{A}^{-1}(y_{h+1})\\
\leq&\frac{\mathfrak{A}^{-1}(\beta)}{\beta^{\frac{1}{n}}}c2^{h}RH^{-1}\omega_{n}^{-1}\mathfrak{A}^{-1}(c\gamma\mathfrak{\widehat{A}}(R^{-1}H)\mathfrak{A}_{*}(2^{h+2})R^ny_h)
\cdot c_3y_{h}\mathfrak{A}^{-1}\bigg(\frac{y_{h+1}}{c_2R^ny_{h}}\bigg)\\
\leq&\frac{\mathfrak{A}^{-1}(\beta)}{\beta^{\frac{1}{n}}}c2^{h}RH^{-1}\omega_{n}^{-1}\mathfrak{A}^{-1}(c\gamma\mathfrak{\widehat{A}}(R^{-1}H)
\mathfrak{A}_{*}(2^{h+2})R^ny_h)\cdot c_4y_{h}\mathfrak{A}^{-1}(R^{-n})\\
\leq&\frac{\mathfrak{A}^{-1}(\beta)}{\beta^{\frac{1}{n}}}c_52^{h}\omega_{n}^{-1}\mathfrak{A}^{-1}(c\gamma
\mathfrak{A}_{*}(2^{h+2})y_h)\cdot y_{h}\mathfrak{A}^{-1}(R^{n})\mathfrak{A}^{-1}(R^{-n})\\
&\quad\quad\cdot(RH^{-1})\mathfrak{A}^{-1}(\mathfrak{\widehat{A}}(R^{-1}H))\\
\leq&\frac{\mathfrak{A}^{-1}(\beta)}{\beta^{\frac{1}{n}}}c_52^{h}\omega_{n}^{-1}\mathfrak{A}^{-1}(c\gamma
\mathfrak{A}_{*}(2^{h+2})y_h)\cdot y_{h}\cdot M_0^2,
\end{aligned}
\end{equation*}
and therefore
\begin{equation*}
\begin{aligned}
\mathfrak{A}_*^{-1}(\beta y_{h+1})
\leq&\frac{\mathfrak{A}^{-1}(\beta)}{\beta^{\frac{1}{n}}}c_62^{h}\omega_{n}^{-1}\mathfrak{A}^{-1}(c\gamma\mathfrak{A}_{*}(2^{h+2})y_h)\cdot\frac{y_{h}}{y_{h+1}}\\
\leq&\frac{\mathfrak{A}^{-1}(\beta)}{\beta^{\frac{1}{n}}}c2^{h}\mathfrak{A}^{-1}(c\mathfrak{A}_{*}(2^{h+2})y_h).
\end{aligned}
\end{equation*}

Then from the above inequality we conclude
\begin{equation}
y_{h+1}\leq\frac{1}{\beta}\mathfrak{A}_*\bigg(\frac{\mathfrak{A}^{-1}(\beta)}{\beta^{\frac{1}{n}}}
c2^{h}\mathfrak{A}^{-1}(c\mathfrak{A}_*(2^{h+2})y_h)\bigg),\,\forall\,\beta>0.
\end{equation}
Choose $\theta\leq\min\{y_0^*,\frac{1}{2}\cdot4^{-n}\omega_{n}\}$. Then \eqref{key1.1} yields that $y_0\leq\theta\leq y_0^*$. By Lemma \ref{sequence} we can get $y_h\downarrow0=\varepsilon_0$ as $h\rightarrow\infty$, which contradicts to $\varepsilon_0>0$.

The conclusion of the \textbf{Claim} means that \eqref{key1.3} holds. The proof is completed.
\end{proof}

\vspace{0.3cm}

\begin{lem}\label{keylemma2}
Let $A\in N(\Omega)\cap\mathscr{A}$ satisfy $(A_1)$ and $(C_1)$. For any given $\theta>0$ there is a natural number $s=s(\theta,n,A,\gamma)>2$ such that either
\begin{equation}\label{2.1}
\psi\leq2^{s}\mathfrak{\widehat{A}}^{-1}\bigg(\frac{\gamma_{1}}{\gamma}\bigg)R
\end{equation}
or
\begin{equation}\label{2.2}
|\Omega_{k^{0},\frac{R}{2}}|\leq\theta R^{n},
\end{equation}
holds, where $k^{0}=\max_{B_{R}}w(x)-2^{-s+1}\psi$.
\end{lem}

\vspace{0.3cm}

\begin{proof}
If equation \eqref{2.1} does not hold, i.e.
\begin{equation}\label{2.1.1}
\psi>2^{s}\mathfrak{\widehat{A}}^{-1}\bigg(\frac{\gamma_{1}}{\gamma}\bigg)R
\end{equation}
where the nature number s will be determined later.\\

For $t=0,1,2\ldots,s-1$, set \\
\begin{equation*}k_{t}=\max_{B_{R}}w(x)-2^{-t}\psi,\quad D_{t}=\Omega_{k_{t},\frac{R}{2}}\setminus\Omega_{k_{t+1},\frac{R}{2}}.
\end{equation*}

Applying inequality \eqref{keyeq} with $\rho=R,\rho-\sigma\rho=\frac{R}{2},k=k_{t}$ for $t=0,1,2,\ldots,s-2$, and by \eqref{2.1.1} we obtain
\begin{equation}\label{key5}
\begin{aligned}[b]
&\int_{D_{t}}A(x,|\nabla w|)\,\mathrm{d}x\\
\leq& \int_{\Omega_{k_{t},\frac{R}{2}}}A(x,|\nabla w|)\,\mathrm{d}x\\
\leq& \gamma\int_{\Omega_{k_{t},R}}A\bigg(x,\frac{w-k_{t}}{\frac{R}{2}}\bigg)\,\mathrm{d}x+\gamma_{1}|\Omega_{k_{t},R}|\\
\leq& \gamma\int_{\Omega_{k_{t},R}}A\bigg(x,\frac{2}{R}2^{-t}\psi\bigg)\,\mathrm{d}x+\gamma_{1}|\Omega_{k_{t},R}|\\
\leq& C_{1}\gamma\mathfrak{\widehat{A}}(R^{-1}2^{-t}\psi)|\Omega_{k_{t},R}|+\gamma_{1}|\Omega_{k_{t},R}|\\
\leq& c\gamma\mathfrak{\widehat{A}}(R^{-1}2^{-t}\psi)R^{n}\quad(t=0,1,\cdots ,s-2).
\end{aligned}
\end{equation}

Applying inequality \eqref{lem3.1} to $l=k_{t+1},k=k_{t}$ and $\rho=\frac{R}{2}$ we get
\begin{equation}\label{key6}
\begin{aligned}[b]
&|\Omega_{k_{s-1},\frac{R}{2}}|^{1-\frac{1}{n}}\\
\leq&|\Omega_{k_{t+1},\frac{R}{2}}|^{1-\frac{1}{n}}\\
\leq&\frac{C}{(k_{t+1}-k_{t})}\frac{(\frac{R}{2})^{n}}{|B_{\frac{R}{2}}-\Omega_{k_{t},\frac{R}{2}}|}\int_{D_{t}}|\nabla w|\,\mathrm{d}x\\
\leq& C(2^{-(t+1)}\psi)^{-1}\bigg(\frac{R}{2}\bigg)^{n}\frac{1}{|B_{\frac{R}{2}}-\Omega_{k_{t},\frac{R}{2}}|}2|\nabla w|_{A;D_{t}}\cdot|1|_{\widetilde{A};D_{t}}\\
\leq&C_12^{t}\psi^{-1}R^{n}(2^{-1}\omega_{n}(2^{-1}R)^{n})^{-1}\mathfrak{A}^{-1}\bigg(\int_{D_{t}}A(x,|\nabla w|)\,\mathrm{d}x\bigg)\cdot\widehat{\widetilde{\mathfrak{A}}}^{-1}\bigg(\int_{D_{t}}\widetilde{A}(x,1)\,\mathrm{d}x\bigg)\\
\leq&c2^{t}\psi^{-1}\mathfrak{A}^{-1}\bigg(\int_{D_{t}}A(x,|\nabla w|)\,\mathrm{d}x\bigg)\cdot\widehat{\widetilde{\mathfrak{A}}}^{-1}(|D_{t}|).
\end{aligned}
\end{equation}

From Lemma \ref{compare}, Remark \ref{rem1}, $\mathfrak{A}^{-1},\mathfrak{A}\in\Delta_{\mathbb{R}^+}$, \eqref{key6} and \eqref{key5} we can conclude
\begin{equation}\label{key7.1}
\begin{aligned}
|\Omega_{k_{s-1},\frac{R}{2}}|^{1-\frac{1}{n}}
\leq&c2^{t}\psi^{-1}\mathfrak{A}^{-1}(c\gamma\mathfrak{\widehat{A}}(R^{-1}2^{-t}\psi)R^{n})\cdot\widehat{\widetilde{\mathfrak{A}}}^{-1}(|D_{t}|)\\
\leq&c(R2^{t}\psi^{-1})\mathfrak{A}^{-1}(\mathfrak{\widehat{A}}(R^{-1}2^{-t}\psi))R^{-1}\mathfrak{A}^{-1}(c\gamma R^{n})\cdot\widehat{\widetilde{\mathfrak{A}}}^{-1}(|D_{t}|)\\
\leq& c_{1}R^{-1}\mathfrak{A}^{-1}(c\gamma R^{n})\cdot\widehat{\widetilde{\mathfrak{A}}}^{-1}(|D_{t}|).
\end{aligned}
\end{equation}
From \eqref{key7.1} and $\widehat{\widetilde{\mathfrak{A}}}\in\Delta_{R^{+}}$  we have
\begin{equation}\label{key7}
\widehat{\widetilde{\mathfrak{A}}}(|\Omega_{k_{s-1},\frac{R}{2}}|^{1-\frac{1}{n}})\leq c_0\widehat{\widetilde{\mathfrak{A}}}(c_{1}R^{-1}\mathfrak{A}^{-1}(c\gamma R^{n}))\cdot|D_{t}|.
\end{equation}
Summing \eqref{key7} with respect to $t$ from 0 to $s-2$ and noting that
\begin{equation*}
\sum_{t=0}^{s-2}|D_{t}|\leq|B_{\frac{R}{2}}|=\omega_{n}\bigg(\frac{R}{2}\bigg)^n,
\end{equation*}
we obtain
 \begin{equation}\label{key8}
 \begin{aligned}
 \widehat{\widetilde{\mathfrak{A}}}(|\Omega_{k_{s-1},\frac{R}{2}}|^{1-\frac{1}{n}})
\leq \frac{c_{0}2^{-n}\omega_n}{s-1}\widehat{\widetilde{\mathfrak{A}}}(c_{1}R^{-1}\mathfrak{A}^{-1}(c\gamma R^{n}))R^n.
\end{aligned}
 \end{equation}
By Lemma \ref{compare}, $\mathfrak{A}^{-1},\mathfrak{A}\in\Delta_{\mathbb{R}^+}$, \eqref{key8} and Remark \ref{rem1}, we can get
 \begin{equation*}
 \begin{aligned}
 |\Omega_{k_{s-1},\frac{R}{2}}|^{1-\frac{1}{n}}
 \leq& M_1\widehat{\widetilde{\mathfrak{A}}}^{-1}\bigg(\frac{c_{0}2^{-n}\omega_n}{s-1}\bigg)c_{1}R^{-1}\mathfrak{A}^{-1}(c\gamma R^{n})\widehat{\widetilde{\mathfrak{A}}}^{-1}(R^n)\\\
 =& M_{1}\widehat{\widetilde{\mathfrak{A}}}^{-1}\bigg(\frac{c_{0}2^{-n}\omega_n}{s-1}\bigg)c_{1}R^{-1}\mathfrak{A}^{-1}(c\gamma R^{n})R^n\mathfrak{A}^{-1}(R^{-n})\\
 \leq&M_{2}\widehat{\widetilde{\mathfrak{A}}}^{-1}\bigg(\frac{c_{0}2^{-n}\omega_n}{s-1}\bigg)c_{1}R^{n-1}\mathfrak{A}^{-1}(c\gamma)
 \mathfrak{A}^{-1}(R^{n})\mathfrak{A}^{-1}(R^{-n})\\
 \leq& M_{3}\widehat{\widetilde{\mathfrak{A}}}^{-1}\bigg(\frac{c_{0}2^{-n}\omega_n}{s-1}\bigg)R^{n-1},
 \end{aligned}
 \end{equation*}
 and therefore
 \begin{equation}\label{key9}
 |\Omega_{k_{s-1},\frac{R}{2}}|\leq\bigg( M_{3}\widehat{\widetilde{\mathfrak{A}}}^{-1}\bigg(\frac{c_{0}2^{-n}\omega_n}{s-1}\bigg)\bigg)^{\frac{n}{n-1}}R^{n}.
 \end{equation}
where $M_3=M_3(\gamma,n,A)$ and $c_0=c_0(A)$ are constants.

 Now we choose a natural number $s$ such that $s-1>c_{0}2^{-n}\omega_n$ and
\begin{equation*}
\bigg( M_{3}\widehat{\widetilde{\mathfrak{A}}}^{-1}\bigg(\frac{c_{0}2^{-n}\omega_n}{s-1}\bigg)\bigg)^{\frac{n}{n-1}}<\theta,
\end{equation*}
 note that $s=s(\theta,n,M_3,c_0,\widehat{\widetilde{\mathfrak{A}}}^{-1})=s(n,\gamma,A,\theta)$. For such $s$, if \eqref{2.1} does not hold, then from \eqref{key9} we get
\begin{equation*}
|\Omega_{k_{s-1},\frac{R}{2}}|\leq \theta R^{n},
\end{equation*}
 which means $|\Omega_{k^{0},\frac{R}{2}}|\leq\theta R^{n}$, i.e. \eqref{2.2} holds since
\begin{equation*}
k^{0}=\max_{B_{R}}w(x)-2^{-s+1}\psi=k_{s-1}.
\end{equation*}
\end{proof}

\vspace{0.3cm}

\begin{lem}\label{keylemma3}
There is a number $s=s(n,A,\gamma)>2$ such that
\begin{equation}\label{3.1}
\psi\leq 2^{s}\max\bigg\{\max_{B_{R}}w(x)-\max_{B_{\frac{R}{4}}}w(x),R \mathfrak{\widehat{A}}^{-1}\bigg(\frac{\gamma_{1}}{\gamma}\bigg)\bigg\}
\end{equation}
for any $R \in (0,R_{0}]$.
\end{lem}

\vspace{0.3cm}

\begin{proof}
Let $R\in (0,R_{0}]$ and let $\theta=\theta(n,A,\gamma)$ be the constant as in Lemma \ref{keylemma1}. Applying Lemma \ref{keylemma2} to this $\theta$ we can find a constant $s=s(\theta,n,A,\gamma)=s(n,A,\gamma)$ such that at least one of \eqref{2.1} and \eqref{2.2} holds.

If \eqref{2.1} holds, then \eqref{3.1} is obviously true.

Now assume that \eqref{2.1} does not hold. Then Lemma \ref{keylemma2}, \eqref{2.2} holds, i.e.
\begin{equation*}
|\Omega_{k^{0},\frac{R}{2}}|\leq\theta R^{n},
\end{equation*}
where $k^{0}=\max_{B_{R}}w(x)-2^{-s+1}\psi$.
Set $H=2^{-s+1}\psi$. Since \eqref{2.1} does not hold, we get
\begin{equation*}
H>2^{-s+1}\cdot2^{s}R\mathfrak{\widehat{A}}^{-1}\bigg(\frac{\gamma_{1}}{\gamma}\bigg)>R\mathfrak{\widehat{A}}^{-1}\bigg(\frac{\gamma_{1}}{\gamma}\bigg).
\end{equation*}
By Lemma \ref{keylemma1}
\begin{equation*}
\max_{B_{R/4}}w(x)\leq k^{0}+\frac{H}{2}=\max_{B_{R}}w(x)-2^{-s}\psi,
\end{equation*}
i.e.
\begin{equation*}
\psi\leq2^{s}\bigg(\max_{B_{R}}w(x)-\max_{R/4}w(x)\bigg).
\end{equation*}
This shows that equation \eqref{3.1} holds.

Lemma \ref{keylemma3} is proved.
\end{proof}

\vspace{0.3cm}

\begin{lem}\label{keylemma4}
For any $R\in(0,R_{0}]$ at least one of the following two inequalities holds
\begin{equation}\label{4.1}
\text{osc}\{u,B_{R}\}\leq\tau2^{s}\mathfrak{\widehat{A}}^{-1}\bigg(\frac{\gamma_{1}}{\gamma}\bigg)R,
\end{equation}
\begin{equation}\label{4.2}
\text{osc}\{u,B_{R/4}\}\leq(1-\tau^{-1}2^{-s})\text{osc}\{u,B_{R}\},
\end{equation}
where $\tau=\max\{2,\frac{2}{\delta}\}$, and $s=s(n,A,\gamma)$ is the constant as in Lemma \ref{keylemma3}.
\end{lem}

\vspace{0.3cm}

\begin{proof}
By Lemma \ref{keylemma3}, equation \eqref{3.1} holds. Then at least one of the following two inequalities holds
\begin{equation}\label{4.3}
\psi\leq 2^{s}R\mathfrak{\widehat{A}}^{-1}\bigg(\frac{\gamma_{1}}{\gamma}\bigg),
\end{equation}
\begin{equation}\label{4.4}
\psi\leq 2^{s}\bigg(\max_{B_{R}}w(x)-\max_{R/4}w(x)\bigg).
\end{equation}

When equation \eqref{4.3} holds we have
\begin{equation*}
\text{osc}\{u;B_{R}\}=\tau\psi\leq\tau2^{s}R\mathfrak{\widehat{A}}^{-1}\bigg(\frac{\gamma_{1}}{\gamma}\bigg),
\end{equation*}
i.e. \eqref{4.1} holds.

When \eqref{4.4} holds we have
\begin{equation*}
\begin{aligned}[b]
\text{osc}\{u;B_{R}\}=&\tau\psi\\
\leq&\tau2^{s}\bigg(\max_{B_{R}}w(x)-\max_{R/4}w(x)\bigg)\\
\leq&\tau2^{s}\bigg(\max_{B_{R}}w(x)-\max_{R/4}w(x)-\min_{B_{R}}w(x)+\min_{R/4}w(x)\bigg)\\
\leq&\tau2^{s}[\text{osc}\{w;B_{R}\}-\text{osc}\{w;B_{R/4}\}].\\
=&\tau2^{s}[\text{osc}\{u;B_{R}\}-\text{osc}\{u;B_{R/4}\}],
\end{aligned}
\end{equation*}
which implies that \eqref{4.2} holds.

Lemma \ref{keylemma4} is proved.
\end{proof}

\vspace{0.3cm}

\textbf{The proof of Theorem \ref{key thm}.} By Lemma \ref{leylemma} and it follows from Lemma \ref{keylemma4} that $u\in C^{0,\alpha}(B_{R_{0}}(x_{0}))$ where
\begin{equation*}
\alpha=\min\{1,-\log_{4}(1-\tau^{-1}2^{-s})\}=-\log_{4}(1-\tau^{-1}2^{-s})=\alpha(n,A,\gamma,\delta).
\end{equation*}
By the arbitrarily of $x_{0}\in\Omega$ we have $u\in C^{0,\alpha}(\Omega)$. The proof of Theorem \ref{key thm} is completed.

\vspace{0.3cm}

\section{Application to minimizers}\label{Sec4}

Consider the integral functionals as follows
\begin{equation}\label{appf}
E(v)=E(v,\Omega)=\int_\Omega f(x,v(x),\nabla
v(x))\,\mathrm{d}x,
\end{equation}
where $v\in W^{1,A}(\Omega)$ and $f(x,s,z)$ is a Carath\'eodory
function on $\Omega\times\mathbb{R}\times\mathbb{R}^n$ satisfying
\begin{equation}\label{appfC}
A\big(x,\sum_{i=1}^n|z_i|\big)-B(x,s)-b\leq f(x,s,z)\leq
a\bigg(A\big(x,\sum_{i=1}^n|z_i|\big)+B(x,s)+b\bigg)
\end{equation}
with $a$ and $b$ being non-negative constants, $A\in N(\Omega)\cap\mathscr{A}$ satisfying $(A_1)$ and $(C_1)$(see In Section \ref{Sec3}), $N(\Omega)\ni
B\preccurlyeq A_*$ satisfying the following $(B_1)$-$(B_2)$:
\begin{enumerate}
\item[$(B_1)$] There exists a strictly increasing
differentiable function
$\mathfrak{B}:[0,+\infty)\rightarrow[0,+\infty)$ satisfying
\begin{equation*}
B(x,\alpha t)\geq \mathfrak{B}(\alpha)B(x,t), \,\forall\,
\alpha\geq0,t\in\mathbb{R},x\in\Omega;
\end{equation*}
\item[$(B_2)$] There exists a constant
$T_{B,\Omega}>0$ such that $B(x,T_{B,\Omega})\geq1$ for any
$x\in\overline\Omega.$
\end{enumerate}

\vspace{0.3cm}

\begin{defi}
A function $u\in W_{loc}^{1,A}(\Omega)$ is said to be a local minimizer of $E$ if
\begin{equation}\label{appmin}
E(u;\text{supp}\varphi)\leq E(u+\varphi;\text{supp}\varphi)\text{ for any }\varphi\in W_0^{1,A}(\Omega)\text{ with } \text{supp}\varphi\subset\subset\Omega.
\end{equation}
\end{defi}

\vspace{0.3cm}

\begin{thm}\label{appthma}
Let $f$ satisfy the growth condition
\eqref{appfC}. If $u\in W^{1,A}(\Omega)$ is a
local minimizer for the functional \eqref{appf},
then
\begin{enumerate}
\item[(i)]$u\in L_{\text{loc}}^\infty(\Omega)$;
\item[(ii)]$u\in C^{0,\alpha}(\Omega)$ in which the constant $\alpha\in(0,1)$ depends only  on $n$, $A$ and $a$.
\end{enumerate}
\end{thm}

\vspace{0.3cm}

The assertion (i) of Theorem \ref{appthma} has been obtained by Theorem 4.1 of \cite{liuandyao}.

The assertion (ii) of Theorem \ref{appthma} follows from Theorem \ref{key thm} and the following Lemma \ref{apple}.

\vspace{0.3cm}

\begin{lem}\label{apple}
Under the assumptions of Theorem \ref{appthma}, for any open set $\Omega_{1}$ in $\Omega$ with $\overline{\Omega}_{1}\subset\Omega$, there exist positive constants $M$, $\gamma$ and $\gamma_{1}$ such that $u\in\mathscr{B}(\Omega_{1},M,\gamma,\gamma_{1},1)$ in which the constant $\gamma=\gamma(a,A)$ is independent of $\Omega_1$.
\end{lem}

\vspace{0.3cm}

\begin{proof}
Let $\Omega_{1}$ be an open set in $\Omega$ with $\Omega_{1}\subset\Omega$. By assertion $(i)$ of Theorem \ref{appthma} there is a constant $M>0$ such that
\begin{equation}\label{applea}
\max_{\overline{\Omega}_{1}}u(x)\leq M.
\end{equation}

From \eqref{appfC} it follow that there is a positive constant $c=c(a,b,B_{+},M)=c(a,b,A,M)$ such that
\begin{equation}\label{appleb}
A\big(x,\sum_{i=1}^n|z_i|\big)-c\leq f(x,u,z)\leq
aA\big(x,\sum_{i=1}^n|z_i|\big)+c,\,\forall\,x\in\overline{\Omega}_{1},\,\ |u|\leq4M,\,\ z\in R^{n},
\end{equation}
where $B_{+}=\max_{\overline{\Omega}_{1}}B(x,M)\leq\max_{\overline{\Omega}_{1}}{A_*(x,M)+C}\leq C_0$.

Let $B_{s}\subset\Omega_{1}$, $0<t<s$, $k\geq-2M$, $w(x)=\max\{u(x)-k,0\}$. Choose $\eta\in C^{\infty}$ with $\text{supp}\eta\subset B_s$, $0\leq\eta\leq1$, $\eta\equiv1$\,\ on $B_t,\,\ |D\eta|\leq2(s-t)^{-1}$. Set $v=u-\eta w$. Then for $x\in\overline{\Omega}_{1}$
\begin{equation*}
|v(x)|\leq|u(x)|+|w(x)|\leq M+3M=4M.
\end{equation*}
By the minimality of $u$ and \eqref{appleb} we get
\begin{equation}\label{applec}
\begin{aligned}
&\int_{\Omega_{k,s}}A(x,|\nabla u|)\,\mathrm{d}x-c|\Omega_{k,s}|\\
\leq& a\int_{\Omega_{k,s}}A(x,|\nabla v|)\,\mathrm{d}x+c|\Omega_{k,s}|\\
\leq& \frac{a}{2}\mathfrak{\widehat{A}}(2)\bigg[\int_{\Omega_{k,s}}A(x,|\nabla u|(1-\eta))\,\mathrm{d}x+\int_{\Omega_{k,s}}A(x,|\nabla\eta|(u-k))\,\mathrm{d}x\bigg]+c|\Omega_{k,s}|\\
\leq& \frac{a}{2}\mathfrak{\widehat{A}}(2)\int_{\Omega_{k,s}\setminus\Omega_{k,t}}A(x,|\nabla u|)\,\mathrm{d}x+\frac{a}{2}(\mathfrak{\widehat{A}}(2))^2\int_{\Omega_{k,s}}A\bigg(x,\frac{u-k}{s-t}\bigg)\,\mathrm{d}x+c|\Omega_{k,s}|.
\end{aligned}
\end{equation}
Therefore
\begin{equation}\label{appled}
\begin{aligned}
&\int_{\Omega_{k,t}}A(x,|\nabla u|)\,\mathrm{d}x\\
\leq& \int_{\Omega_{k,s}}A(x,|\nabla u|)\,\mathrm{d}x\\
\leq&  \frac{a}{2}\mathfrak{\widehat{A}}(2)\int_{\Omega_{k,s}\setminus\Omega_{k,t}}A(x,|\nabla u|)\,\mathrm{d}x+ \frac{a}{2}(\mathfrak{\widehat{A}}(2))^2\int_{\Omega_{k,s}}A\bigg(x,\frac{u-k}{s-t}\bigg)\,\mathrm{d}x+2c|\Omega_{k,s}|.
\end{aligned}
\end{equation}
Adding $\frac{a}{2}\mathfrak{\widehat{A}}(2)\int_{\Omega_{k,t}}A(x,|\nabla u|)\,\mathrm{d}x$ to both sides of \eqref{appled} we get
\begin{equation}\label{applee}
\begin{aligned}
&\int_{\Omega_{k,t}}A(x,|\nabla u|)\,\mathrm{d}x\\
\leq& \theta\int_{\Omega_{k,s}}A(x,|\nabla u|)\,\mathrm{d}x+c_{1}\int_{\Omega_{k,s}}A\bigg(x,\frac{u-k}{s-t}\bigg)\,\mathrm{d}x+c_{2}|\Omega_{k,s}|,
\end{aligned}
\end{equation}
where $\theta=\frac{a\mathfrak{\widehat{A}}(2)}{2+a\mathfrak{\widehat{A}}(2)}<1$, $c_{1}=c_{1}(a,A)$ and $c_{2}=c_{2}(c,a,A)=c_{2}(a,b,A,M)$ are positive constants.
Using the similar method that was used in the proof of Lemma 4.1 of \cite{liuandyao}, from  \eqref{applee} we can deduce that
\begin{equation}\label{applef}
\int_{\Omega_{k,\sigma R}}A(x,|\nabla u|)\,\mathrm{d}x\leq
\gamma\int_{\Omega_{k,R}}A\bigg(x,\frac{u-k}{(1-\sigma)R}\bigg)\,\mathrm{d}x
+\gamma_{1}|\Omega_{k,s}|,
\end{equation}
for $B_{R}\subset\Omega_{1}$, $\sigma\in(0,1)$ and $k\geq-2M$, where $\gamma=\gamma(a,A)$ and $\gamma_{1}=\gamma_{1}(a,b,A,M)$ are positive constants.

Set $\delta=1$. Then for every $B_{\rho}\subset\Omega_{1}$
\begin{equation*}
\max_{B_{\rho}}u(x)-M\geq-2M,
\end{equation*}
and consequently \eqref{applef} holds for $k\geq\max_{B_{R}}u(x)-\delta M$ with $\delta=1$.

It is easy to see that \eqref{applef} holds when $u$ is replaced by $-u$.
So by Definition \ref{keyde}, $u\in\mathscr{B}(\Omega_{1},M,\gamma,\gamma_{1},1)$.  Lemma \ref{apple} is proved.
\end{proof}

\vspace{0.3cm}

\textbf{The proof of Theorem \ref{appthma}.} In the case $A$  satisfies $(A_1)$ and $(C_1)$, by Lemma \ref{apple} and Theorem \ref{key thm}, $u\in C^{0,\alpha}(\overline{\Omega}_1)$ for every $\overline{\Omega}_1\subset\Omega$ and the constant $\alpha\in(0,1)$ is independent of $\Omega_1$. Hence $u\in C^{0,\alpha}(\Omega)$. The proof of Theorem \ref{appthma} is completed.

\vspace{0.3cm}

\section{Application to fully nonlinear elliptic equations}\label{Sec5}

In this section, we consider the local H\"{o}lder continuity of weak solutions of a kind of
fully nonlinear elliptic equation. Since we only consider the local properties of the
weak solutions, without loss of generality, we suppose that $\Omega$ is a bounded smooth
domain in $\mathbb{R}^n$.

\vspace{0.3cm}

Consider the second order fully nonlinear elliptic equation as follows
\begin{equation}\label{appq}
\text{div} L(x,u,\nabla u)+F(x,u,\nabla u)=0,\, \forall\,x\in\Omega,
\end{equation}
where $L:\Omega
\times\mathbb{R}\times\mathbb{R}^n\rightarrow\mathbb{R}^n$, $F:\Omega
\times\mathbb{R}\times\mathbb{R}^n\rightarrow\mathbb{R}^1$, and $u:\Omega\rightarrow\mathbb{R}$.

\vspace{0.3cm}

Suppose equation \eqref{appq} satisfies the following growth conditions:
\begin{equation}\label{appqca}
L(x,u,z)z\geq a_0A(x,|z|)-bB(x,u)-c,
\end{equation}
\begin{equation}\label{appqcb}
|L(x,u,z)|\leq a_1\widetilde A^{-1}A(x,|z|)+b\widetilde A^{-1}B(x,u)+c,
\end{equation}
\begin{equation}\label{appqcc}
|F(x,u,z)|\leq a_2\widetilde B^{-1}A(x,|z|)+b\widetilde B^{-1}B(x,u)+c,
\end{equation}
where $a_0,a_1,a_2,b,c$ are positive constants, $A\in N(\Omega)\cap\mathscr{A}$ satisfies
$(A_1)$ and $(C_1)$, $N(\Omega)\ni
B\preccurlyeq A_*$ satisfies $(B_1^+)$-$(B_2^+)$.
\begin{enumerate}
\item[$(B_1^+)$] There exists a strictly increasing
differentiable function
$\mathfrak{B}:[0,+\infty)\rightarrow[0,+\infty)$ satisfying
\begin{equation}
\alpha\mathfrak{B}'(\alpha)>\mathfrak{B}(\alpha)
\end{equation}
such that
\begin{equation*}
B(x,\alpha t)\geq \mathfrak{B}(\alpha)B(x,t), \,\forall\,
\alpha\geq0,t\in\mathbb{R},x\in\Omega;
\end{equation*}
\item[$(B_2^+)$] There exists a constant
$T_{B,\Omega}>0$ such that $B(x,T_{B,\Omega})\geq1$ for any
$x\in\overline\Omega$.
\end{enumerate}

\vspace{0.3cm}

\begin{defi}
$u\in W^{1,A}(\Omega)$ is said to be a weak solution of \eqref{appq} if
\begin{equation}\label{app w s}
\int_\Omega L(x,u,\nabla u)\nabla v\,\mathrm{d}x-\int_\Omega F(x,u,\nabla u) v\,
\mathrm{d}x=0
\end{equation}
for any $v\in W_0^{1,A}(\Omega)$.
\end{defi}

\vspace{0.3cm}

The local bounded regularity of weak solutions of \eqref{appq} satisfying \eqref{appqca}-\eqref{appqcc} has been obtained by Theorem 5.1 of \cite{liuandyao}. Now we discuss the H\"{o}lder continuity of weak solutions of \eqref{appq}.

\vspace{0.3cm}

\begin{thm}\label{appthmb}
Let the growth conditions \eqref{appqca}-\eqref{appqcc} hold. If $u\in W^{1,A}(\Omega)$ is a
weak solution of \eqref{appq} and $\max_{\Omega}|u(x)|\leq M$, then $u\in\mathscr{B}(\Omega,M,\gamma,\gamma_{1},\delta)$ in which $\gamma=\gamma(a_{0},a_{1},A)$, $\gamma_1=\gamma_{1}(a_{0},a_{2},b,c,A,M)$, $\delta=\min\{\frac{a_0}{4a_{2}M},2\}$.
\end{thm}

\vspace{0.3cm}

\begin{proof}
Let  $u$ be a weak solution of \eqref{appq}. For arbitrary balls $\overline{B}_{s}\subset B_{t}\subset\Omega$ and pick a function $\xi\in C^{\infty}(\Omega)$ such that
\begin{equation}\label{app q th a}
0\leq\xi\leq1,\quad \text{supp}\xi\subset B_{t},\quad\xi\equiv1\text{ on }B_{s},\quad|D\xi|\leq\frac{2}{t-s}.
\end{equation}
Let
\begin{equation}\label{app q th b}
\delta=\min\bigg\{\frac{a_0}{4a_{2}M},2\bigg\},
\end{equation}
and $v=\mathfrak{A}(\xi)\max\{u-k,0\}\in W_{0}^{1,A}(\Omega)$, where
\begin{equation}\label{app q th c}
k\geq\max_{B_{t}}u(x)-\delta M.
\end{equation}

By \eqref{app w s} we obtain
\begin{equation}\label{app q th d}
\begin{aligned}
\int_{\Omega_{k,t}}\mathfrak{A}(\xi) L(x,u,\nabla u)&\cdot\nabla u\,\mathrm{d}x+\int_{\Omega_{k,t}}(u-k)L(x,u,\nabla u)\cdot\nabla \mathfrak{A}(\xi)\,\mathrm{d}x\\
&\quad\quad-\int_{\Omega_{k,t}}\mathfrak{A}(\xi)(u-k)F(x,u,\nabla u)\,\mathrm{d}x=0.
\end{aligned}
\end{equation}
From \eqref{appqca}--\eqref{appqcc} and \eqref{app q th d},  it follows that
\begin{equation}\label{app q th e}
\begin{aligned}
&a_0\int_{\Omega_{k,t}}A(x,|\nabla u|)\mathfrak{A}(\xi)\,\mathrm{d}x\leq b\int_{\Omega_{k,t}}B(x,|u|)\mathfrak{A}(\xi)\,\mathrm{d}x+c\int_{\Omega_{k,t}}\mathfrak{A}(\xi)\,\mathrm{d}x\\
&\quad\quad+a_1\int_{\Omega_{k,t}}\widetilde A^{-1}A(x,|\nabla u|)|\nabla \mathfrak{A}(\xi)|(u-k)\,\mathrm{d}x\\
&\quad\quad+b\int_{\Omega_{k,t}}\widetilde A^{-1}B(x,|u|)|\nabla \mathfrak{A}(\xi)|(u-k)\,\mathrm{d}x+c\int_{\Omega_{k,t}}|\nabla \mathfrak{A}(\xi)|(u-k)\,\mathrm{d}x\\
&\quad\quad+a_2\int_{\Omega_{k,t}}\widetilde B^{-1}A(x,|\nabla u|)\mathfrak{A}(\xi)(u-k)\,\mathrm{d}x\\
&\quad\quad+b\int_{\Omega_{k,t}}\widetilde B^{-1}B(x,|u|)\mathfrak{A}(\xi)(u-k)\,\mathrm{d}x+c\int_{\Omega_{k,t}}\mathfrak{A}(\xi)(u-k)\,\mathrm{d}x.
\end{aligned}
\end{equation}

We will estimate each term of the right-hand side of \eqref{app q th e}. By $B\preccurlyeq A_*$ and Remark \ref{rem} we obtain
\begin{equation}\label{app q th f}
\begin{aligned}
b\int_{\Omega_{k,t}}B(x,|u|)\mathfrak{A}(\xi)\,\mathrm{d}x\leq& b\int_{\Omega_{k,t}}B(x,M)\mathfrak{A}(1)\,\mathrm{d}x\\
\leq& b\int_{\Omega_{k,t}}A_*(x,M)\mathfrak{A}(1)\,\mathrm{d}x+b\int_{\Omega_{k,t}}C\mathfrak{A}(1)\,\mathrm{d}x\\
\leq& C_0|\Omega_{k,t}|,
\end{aligned}
\end{equation}
where $C_0=C_0(b,A,M)$, and
\begin{equation}\label{app q th g}
c\int_{\Omega_{k,t}}\mathfrak{A}(\xi)\,\mathrm{d}x\leq c\int_{\Omega_{k,t}}\mathfrak{A}(1)\,\mathrm{d}x= c\mathfrak{A}(1)|\Omega_{k,t}|.
\end{equation}

By the Young inequality, and taking $\epsilon>0$ such that
\begin{equation*}
a_1\widetilde{\mathfrak{A}}(n\epsilon)=\frac{a_0}{4},
\end{equation*}
we deduce from the assumption $n\mathfrak{A}(\alpha)>\alpha\mathfrak{A}'(\alpha)>\mathfrak{A}(\alpha)$ and Lemma \ref{compare} (iii) that
\begin{equation}
\begin{aligned}
&a_1\int_{\Omega_{k,t}}\widetilde A^{-1}A(x,|\nabla u|)|\nabla \mathfrak{A}(\xi)|(u-k)\,\mathrm{d}x\\
=&a_1\int_{\Omega_{k,t}}\widetilde A^{-1}A(x,|\nabla u|)|\nabla \xi|\mathfrak{A}'(\xi)(u-k)\,\mathrm{d}x\\
\leq& a_1\int_{\Omega_{k,t}}\widetilde A\big(x,\epsilon \widetilde A^{-1}A(x,|\nabla u|)\mathfrak{A}'(\xi)\big)\,\mathrm{d}x
+a_1\int_{\Omega_{k,t}}A(x,\epsilon^{-1}|\nabla \xi|(u-k))\,\mathrm{d}x\\
\leq& a_1\int_{\Omega_{k,t}}\widetilde A\bigg(x,\epsilon \widetilde A^{-1}A(x,|\nabla u|)\frac{n\mathfrak{A}(\xi)}{\xi}\bigg)\,\mathrm{d}x
+a_1\int_{\Omega_{k,t}}A(x,\epsilon^{-1}|\nabla \xi|(u-k))\,\mathrm{d}x\\
\leq&a_1\widetilde{\mathfrak{A}}(n\epsilon)\int_{\Omega_{k,t}}A(x,|\nabla u|)\mathfrak{A}(\xi)\,\mathrm{d}x
+a_1\widehat{\mathfrak{A}}\bigg(\frac{2}{\epsilon}\bigg)\int_{\Omega_{k,t}}A\bigg(x,\frac{u-k}{t-s}\bigg)\,\mathrm{d}x\\
=&\frac{a_0}{4}\int_{\Omega_{k,t}}A(x,|\nabla u|)\mathfrak{A}(\xi)\,\mathrm{d}x
+a_1\widehat{\mathfrak{A}}\bigg(\frac{2}{\epsilon}\bigg)\int_{\Omega_{k,t}}A\bigg(x,\frac{u-k}{t-s}\bigg)\,\mathrm{d}x.
\end{aligned}
\end{equation}
By the Young inequality, \eqref{app q th b} and \eqref{app q th c} we obtain that
\begin{equation}
\begin{aligned}
&a_2\int_{\Omega_{k,t}}\widetilde B^{-1}A(x,|\nabla u|)\mathfrak{A}(\xi)(u-k)\,\mathrm{d}x\\
\leq& a_2\int_{\Omega_{k,t}}A(x,|\nabla u|)\mathfrak{A}(\xi)(u-k)\,\mathrm{d}x+a_2\int_{\Omega_{k,t}}B(x,1)\mathfrak{A}(\xi)(u-k)\,\mathrm{d}x\\
\leq& a_2\int_{\Omega_{k,t}}A(x,|\nabla u|)\mathfrak{A}(\xi)(\max_{B_{t}} u(x)-k)\,\mathrm{d}x+a_2\int_{\Omega_{k,t}}B(x,1)\mathfrak{A}(1)(\max_{B_{t}} u(x)-k)\,\mathrm{d}x\\
\leq& a_2\delta M\int_{\Omega_{k,t}}A(x,|\nabla u|)\mathfrak{A}(\xi)\,\mathrm{d}x+a_2\delta M\mathfrak{A}(1)\bigg(\int_{\Omega_{k,t}}A_*(x,1)\,\mathrm{d}x+\int_{\Omega_{k,t}}C\,\mathrm{d}x\bigg)\\
\leq& \frac{a_0}{4}\int_{\Omega_{k,t}}A(x,|\nabla u|)\mathfrak{A}(\xi)\,\mathrm{d}x+C_1|\Omega_{k,t}|,
\end{aligned}
\end{equation}
where $C_1=C_1(a_0,A).$ Using Young's inequality and taking $\varepsilon_{2}\in(0,1)$ such that
\begin{equation*}
b\widehat{\mathfrak{A}}(2\varepsilon_2)=1,
\end{equation*}
we deduce from the assumption $n\mathfrak{A}(\alpha)>\alpha\mathfrak{A}'(\alpha)>\mathfrak{A}(\alpha)$ and Lemma \ref{compare} (iii) that
\begin{equation}
\begin{aligned}
&b\int_{\Omega_{k,t}}\widetilde A^{-1}B(x,|u|)|\nabla\mathfrak{A}(\xi)|(u-k)\,\mathrm{d}x\\
=&b\int_{\Omega_{k,t}}\widetilde A^{-1}B(x,|u|)|\nabla \xi|\mathfrak{A}'(\xi)(u-k)\,\mathrm{d}x\\
\leq& b\int_{\Omega_{k,t}}\widetilde{A}(x,\varepsilon_2^{-1}\widetilde{A}^{-1}B(x,|u|)\mathfrak{A}'(\xi))\,\mathrm{d}x
+b\int_{\Omega_{k,t}}A(x,\varepsilon_2|\nabla \xi|(u-k))\,\mathrm{d}x\\
\leq& b\int_{\Omega_{k,t}}\widetilde{A}\bigg(x,\varepsilon_2^{-1}\widetilde{A}^{-1}B(x,|u|)\frac{n\mathfrak{A}(\xi)}{\xi}\bigg)\,\mathrm{d}x
+b\widehat{\mathfrak{A}}(2\varepsilon_2)\int_{\Omega_{k,t}}A\bigg(x,\frac{u-k}{t-s}\bigg)\,\mathrm{d}x\\
\leq& b\widetilde{\mathfrak{A}}(\varepsilon_2^{-1}n)\int_{\Omega_{k,t}}B(x,M)\mathfrak{A}(\xi)\,\mathrm{d}x
+b\widehat{\mathfrak{A}}(2\varepsilon_2)\int_{\Omega_{k,t}}A\bigg(x,\frac{u-k}{t-s}\bigg)\,\mathrm{d}x\\
\leq&C_2|\Omega_{k,t}|+\int_{\Omega_{k,t}}A\bigg(x,\frac{u-k}{t-s}\bigg)\,\mathrm{d}x,
\end{aligned}
\end{equation}
where $C_2=C_2(b,A,n,M)$.
Using  Young's inequality we get
\begin{equation}
\begin{aligned}
&c\int_{\Omega_{k,t}}|\nabla\mathfrak{A}(\xi)|(u-k)\,\mathrm{d}x\\
=&c\int_{\Omega_{k,t}}|\nabla \xi|\mathfrak{A}'(\xi)(u-k)\,\mathrm{d}x\\
\leq& 2\int_{\Omega_{k,t}}A\bigg(x,\frac{u-k}{t-s}\bigg)\,\mathrm{d}x+2\int_{\Omega_{k,t}}\widetilde A(x,c\mathfrak{A}'(\xi))\,\mathrm{d}x\\
\leq& 2\int_{\Omega_{k,t}}A\bigg(x,\frac{u-k}{t-s}\bigg)\,\mathrm{d}x+\int_{\Omega_{k,t}}\widetilde A\bigg(x,c\frac{n\mathfrak{A}(\xi)}{\xi}\bigg)\,\mathrm{d}x\\
\leq& 2\int_{\Omega_{k,t}}A\bigg(x,\frac{u-k}{t-s}\bigg)\,\mathrm{d}x+\int_{\Omega_{k,t}}\widetilde A(x,cn)\mathfrak{A}(\xi)\,\mathrm{d}x\\
\leq& 2\int_{\Omega_{k,t}}A\bigg(x,\frac{u-k}{t-s}\bigg)\,\mathrm{d}x+C_3|\Omega_{k,t}|,
\end{aligned}
\end{equation}
where $C_3=C_3(c,n,A)$.
Similarly, we have
\begin{equation}
\begin{aligned}
&b\int_{\Omega_{k,t}}\widetilde B^{-1}B(x,|u|)\mathfrak{A}(\xi)(u-k)\,\mathrm{d}x\\
\leq& b\int_{\Omega_{k,t}}B(x,|u|)\,\mathrm{d}x+b\int_{\Omega_{k,t}}B(x,\mathfrak{A}(\xi)(u-k))\,\mathrm{d}x\\
\leq& b\int_{\Omega_{k,t}}B(x,M)\,\mathrm{d}x+b\int_{\Omega_{k,t}}B(x,\delta M\mathfrak{A}(1))\,\mathrm{d}x\\
\leq& b\int_{\Omega_{k,t}}B(x,M)\,\mathrm{d}x+b\int_{\Omega_{k,t}}B\bigg(x,\frac{a_0}{4a_2}\mathfrak{A}(1)\bigg)\,\mathrm{d}x\\
\leq& C_4|\Omega_{k,t}|,
\end{aligned}
\end{equation}
where $C_4=C_4(b,A,M,a_0,a_2)$.
\begin{equation}\label{app q th h}
\begin{aligned}
c\int_{\Omega_{k,t}}\mathfrak{A}(\xi)(u-k)\,\mathrm{d}x\leq c\int_{\Omega_{k,t}}\mathfrak{A}(1)\delta M\,\mathrm{d}x
\leq \frac{c\mathfrak{A}(1)a_0}{4a_2}|\Omega_{k,t}|.
\end{aligned}
\end{equation}

From \eqref{app q th e}-\eqref{app q th h} we conclude that
\begin{equation*}
\int_{\Omega_{k,t}}A(x,|\nabla u|)\mathfrak{A}(\xi)\,\mathrm{d}x\leq \gamma\int_{\Omega_{k,t}}A\bigg(x,\frac{u-k}{t-s}\bigg)\,\mathrm{d}x
+\gamma_{1}|\Omega_{k,t}|.
\end{equation*}
Therefore
\begin{equation}\label{app q th l}
\mathfrak{A}(1)\int_{\Omega_{k,s}}A(x,|\nabla u|)\,\mathrm{d}x\leq \gamma\int_{\Omega_{k,t}}A\bigg(x,\frac{u-k}{t-s}\bigg)\,\mathrm{d}x
+\gamma_{1}|\Omega_{k,t}|,
\end{equation}
for $\overline{B}_{s}\subset B_t\subset\Omega$ and $k$ satisfying \eqref{app q th c}, where $\gamma=\gamma(a_0,a_1,n,A)$ and $\gamma_1=\gamma_1(a_0,a_2,b,c,A,n,M)$
 are positive constants. \eqref{app q th l} shows that $u\in\mathscr{B}(\Omega,M,\gamma,\gamma_1,\delta)$. Theorem \ref{appthmb} is proved.
\end{proof}

\vspace{0.3cm}

From Theorems \ref{appthmb} and \ref{key thm} we obtain

\vspace{0.3cm}

\begin{thm}
Suppose that the assumptions of Theorems \ref{appthmb} hold. If $u\in W^{1,A}(\Omega)$ is a weak solution of \eqref{appq} and $\max_{\Omega}|u(x)|\leq M$, then $u\in C^{0,\alpha}(\Omega)$ in which $\alpha=\alpha(a_0,a_1,\delta,A,n)=\alpha(a_0,a_1,a_2,M,A,n)\in(0,1)$.
\end{thm}

\renewcommand{\baselinestretch}{0.1}
\bibliographystyle{plain}
\bibliography{Ref}

\end{document}